\documentclass[12pt,36pc, a4paper]{amsart}
\usepackage{macs-symimp}
\allowdisplaybreaks[2]
 \date{February 5,  2002}
\begin{document}

\author[an Huef]{Astrid an Huef} \address{School of Mathematics \\ 
  University of New South Wales \\ Sydney, NSW 2052 \\ Australia}
\email{astrid@maths.unsw.edu.au}

\author[Raeburn]{Iain Raeburn}
\address{School of Mathematics and Physical Sciences \\
University of Newcastle \\
Callaghan, NSW 2308 \\
Australia}
\email{iain@maths.newcastle.edu.au}

\author[Williams]{Dana P. Williams}
\address{Department of Mathematics \\
Dartmouth College \\
Hanover, NH 03755-3551 \\
USA}
\email{dana.williams@dartmouth.edu}

\thanks{This research was supported by grants from the Australian Research Council, the National Science Foundation, Dartmouth College, the University of Denver, and the University of Newcastle.}


\title[Proper Actions]{A Symmetric Imprimitivity Theorem for Commuting
  Proper Actions}

\begin{abstract}
  We initiate a careful study of a generalized symmetric imprimitivity
  theory for commuting proper actions of locally compact groups $H$
  and $K$ on a \cs-algebra.
\end{abstract}

\maketitle

\section{Introduction}

Until recently, the various symmetric imprimitivity theorems in the
literature have all been associated to commuting free and proper actions of two
locally compact groups $H$ and $K$ on the left and right of a locally
compact space $P$. The original theorem, due to Green and Rieffel
\cite{rie:pspm82}, says that the crossed products $C_{0}(P/H)\rtimes
K$ and $C_{0}(K\backslash P)\rtimes H$ are Morita equivalent; the most
powerful generalizations give Morita equivalences between crossed products of
induced $C^*$-algebras \cite{rae:ma88,kas} or crossed products of
$C_0(P)$-algebras \cite{hrw:tams00}.

There two ways to prove a symmetric imprimitivity theorem. The first,
used in \cite{rie:pspm82} and \cite{rae:ma88}, is to build module
actions and inner products on spaces of compactly supported functions
and then complete to get a bimodule over the (complete) crossed
products; for the Green-Rieffel theorem, the bimodule $Z=Z_K^H$ is a
completion of $C_c(P)$. The second, used in \cite{cmw} and \cite{kas},
starts from the one-sided equivalence involving just one group, and
bootstraps up by taking crossed products and tensor products; for the
Green-Rieffel theorem, we start with $X:=Z_{\set e}^{H}$ and
$Y:=Z^{\set e}_{K}$, and the tensor product
\begin{equation*}
(X\rtimes K)\tensor_{C_{0}(P)\rtimes (K\times H)} (Y\rtimes H).
\end{equation*}
implements the equivalence. In applications, the first bimodule $Z$ is
more convenient for direct calculations, and the tensor-product
bimodule is useful when we want to bootstrap results from the
one-sided case. It is known that the two bimodules are in fact isomorphic as
imprimitivity bimodules \cite{hrw:tams00}, and this isomorphism is useful, for example, in
settling questions of amenability \cite{hrw:tams00, huerae:mpcps00}.

Both constructions, though, ultimately use imprimitivity bimodules
constructed from algebras of functions, and the algebraic structure in
these bimodules was found by pretty much \emph{ad hoc} methods.  In
\cite{rie:pm88}, Rieffel proposed an alternative, more systematic
approach for building imprimitivity bimodules, based on abstracting
the concepts of proper and free actions to the noncommutative setting.
He described a family of proper saturated actions of a locally
compact group $G$ on a non-commutative $C^*$-algebra $C$ for which
there is a Morita equivalence between the (reduced) crossed product
$C\rtimes_r G$ and a generalized fixed-point algebra $C^G$ of $C$; this
one-sided equivalence is implemented by a bimodule which is
constructed by completing a dense subalgebra of $C$ in a very
particular way.

Pask and Raeburn have recently proved a symmetric imprimitivity
theorem for commuting actions on the Cuntz-Krieger algebras of
directed graphs \cite[Theorem~2.1]{PR}, which appears to be quite
independent of the machinery developed in
\cite{rie:pspm82,rae:ma88,hrw:tams00}. However, the actions considered
in \cite{PR} are proper in Rieffel's sense, and \cite[Theorem~2.1]{PR}
can be formulated as a Morita equivalence of crossed products of
generalized fixed-point algebras. It is therefore tempting to look for
a symmetric imprimitivity theorem for commuting proper actions on a
$C^*$-algebra, and the purpose of the present paper is to formulate and
prove such a theorem. Thus we consider commuting actions $\tau:H\to
\Aut C$ and $\sigma:K\to \Aut C$ of two groups on the same
$C^*$-algebra $C$, and aim to prove that if both actions are proper
and saturated in Rieffel's sense, then we have a Morita equivalence
between the crossed products $C^\tau\rtimes_{\sigma,r} K$ and $C^\sigma\rtimes_{\tau,r} H$. We
want a complete theory: we want a tensor-product bimodule which is
good for bootstrapping arguments, a bimodule which is obtained by
completing a dense subalgebra of $C$, and
an isomorphism between these bimodules. If this new symmetric
imprimitivity theorem requires extra hypotheses, we want to know that
the hypotheses are satisfied in the key examples.

The first step is relatively straightforward. Under some mild
continuity hypotheses which ensure that the various crossed products
make sense, we can start with two applications of Rieffel's theorem
from \cite{rie:pm88} and use the usual bootstrap arguments to obtain a
tensor-product bimodule (Proposition~\ref{prop-tensor-ib}). The second
and third steps are achieved in \S\ref{concrete-result} using the
results of our earlier paper \cite{hrw:xx}. We show that the natural
action of $K$ on the bimodule implementing the equivalence between $C^\tau$
and $C\rtimes_{\tau, r} H$ is proper and saturated in the sense of
\cite{hrw:xx}, and identify the generalized fixed-point algebra
$(C\rtimes_{\tau, r} H)^K$ with $C^\sigma\rtimes_{\tau,r} H$, so that the main
theorems of \cite{hrw:xx} give the desired Morita equivalence between
$C^\tau\rtimes_{\sigma,r} K$ and $C^\sigma\rtimes_{\tau,r} H$ and the isomorphism with the
tensor-product bimodule (Theorem~\ref{thm-symmetric} and Corollary~\ref{cor-ibm}).  The proof of Theorem~\ref{thm-symmetric} raises  substantial technical problems
involving vector-valued integrals whose treatment we defer to two appendices.

Theorem~\ref{thm-symmetric} requires substantial hypotheses 
of the sort needed by Rieffel in \cite{rie:pm88}.
In the final section, we show that the hypotheses in
\S\ref{concrete-result} are often automatically satisfied. More
specifically, we show that if there is an underlying free and proper
space ${}_KP_H$ such that $C_0(P)$ maps bi-equivariantly into $M(C)$,
then Theorem~\ref{thm-symmetric} applies.  Although it is a little
against the spirit of Rieffel's theory to assume the
existence of an underlying proper action on a space, it is a fact that
in all main examples of proper actions of $G$ there is such an
underlying space ${}_GP$.  We discuss this in our concluding Remark~\ref{rem-coactions}, and also  speculate on
possible implications for nonabelian duality.

\section{The Tensor-Product Imprimitivity Bimodule}
\label{tensor-result}

Let $\tau:H\to \Aut C$ and $\sigma:K\to \Aut C$ be commuting actions
of locally compact groups on a \cs-algebra $C$.  We assume that both
$\tau$ and $\sigma$ are proper and saturated with respect to the same
dense invariant $*$-subalgebra $C_0$ of $C$ in the sense of
\cite{{rie:pm88}}.

Applying \cite[Corollary 1.7]{rie:pm88} to $\tau$ gives a
$C\rtimes_{\tau,r}H\sme C^\tau$ imprimitivity bimodule
$\overline{C_0}$, where $C^\tau$ denotes the generalized fixed-point
algebra.  Recall that, by definition, $C^\tau$ is the closure of
\[
D_0:=\sp\set{ \rip C^\tau<b, c>: b,c\in C_0}\subset M(C)^\tau,
\] 
where each $\rip C^\tau<b, c>$ is a uniquely determined element of
$M(C)^\tau$ such that for every $a\in C_0$
\[\int_H a\tau_s(b^*c)\, ds=a \rip C^\tau<b, c>.
\]
Since $\tau$ is saturated,
\[
E_0:=\sp\set{ s\mapsto\Delta_H(s)^{-1/2}b\tau_s(c^*): b,c\in C_0}
\]
is dense in $C\rtimes_{\tau,r}H$.

Throughout we will denote by $X$ a module isomorphic to the dual of
$\metrip{C\rtimes_{\tau,r}H}{(\overline{C_0})}{C^\tau}$, so that $X$
is a $C^\tau\sme (C\rtimes_{\tau,r}H)$-imprimitivity bimodule.
Formally, $X$ is obtained by completing $X_0:=C_0$, where $X_0$ is the
left $D_0$-module with $d\cdot x:=dx$ and $\lip C^\tau<x,y>=\rip
C^\tau< x^*, y^*>$; one can easily check that the map
$\phi:\flat(c)\mapsto c^*$ is an isomorphism of the dual of $\metrip
{C^\tau}{(\overline{C_0})}{}$ onto $\metrip {C^\tau}{X_0}{}$.  We use
the same isomorphism $\phi$ to work out what the formula for the
$\cctrh$-valued inner product on $X_0$ should be:
\begin{align*}
  \rip\cctrh<x,y>(s)&=\brip\cctrh<(x^*)^*,(y^*)^*>(s)\\
  &=\brip\cctrh<\phi(\flat(x^*)),\phi(\flat(y^*))>(s)\\
  &=\brip \cctrh<\flat(x^*),\flat(y^*)>(s) \\
  &=\blip \cctrh<x^*, y^*>(s)= \Delta_H(s)^{-1/2}x^*\tau_s(y).
\end{align*}
Now $X_0=C_0$ completes to give a Morita equivalence between $C^\tau$
and $C\rtimes_{\tau,r}H$, and for $x,y,z\in X_0, d\in D_0$ and $e\in
E_0\subset L^1(H,C)$ the actions and inner products are given by
\begin{align}
  &d\cdot x=dx\text{\ is multiplication in $M(C)$}\label{dual-la}\\
  &x\cdot e=\int_H \tau_s^{-1}(xe(s))\Delta_H(s)^{-1/2}\,ds\label{dual-ra}\\
  &\lip {C^\tau}<x,v>\text{\ is characterized by\ }
  \lip {C^\tau}<x,v>\cdot z=\int_H\tau_s(xv^*)z\, ds\label{dual-lip}\\
  &\rip \cctrh<x, v>(s)=\Delta_H(s)^{-1/2}
  x^*\tau_s(v)\label{dual-rip}.
\end{align}

The first step to obtaining the tensor-product version of the
symmetric imprimitivity theorem is to show that the natural extension
$\bar\sigma$ of $\sigma$ to $M(C)$ leaves $C^\tau$ invariant: if
$x,v,w\in X_0$ then, using \eqref{dual-lip}, we have
\begin{align*}
  \blip {C^\tau}<\sigma_t(x),\sigma_t(v)>\cdot w &=\int_H
  \tau_s(\sigma_t(x)\sigma_t(v)^*)w\, ds
  =\int_H \sigma_t(\tau_s(xv^*))w\, ds\\
  &=\sigma_t\Big(\int_H \tau_s(xv^*)\sigma_t^{-1}(w)\, ds\Big)
  =\sigma_t\big(\lip {C^\tau}<x, v>\cdot\sigma_t^{-1}(w)\big)\\
  &=\bar\sigma_t\big(\lip {C^\tau}< x, v>\big)\cdot w.
\end{align*}
Similarly, we use \eqref{dual-rip} to show that
\begin{align*}
  \brip C\rtimes_{\tau,r}H< \sigma_t(x),\sigma_t(v)>(s) &=
  \Delta_H(s)^{-1/2}\sigma_t(x)^*\tau_s(\sigma_t(v))=
  \Delta_H(s)^{-1/2}\sigma_t(x^*\tau_s(v))\\
  &=\sigma_t(\rip C\rtimes_{\tau,r}H< x,v>(s)) =(\sigma\rtimes\id)_t
  (\rip C\rtimes_{\tau,r}H<x,v>)(s).
\end{align*}
These two calculations show that, algebraically at least,
$(\bar\sigma, \sigma,\sigma\rtimes\id)$ is a potential candidate for
an action of $K$ on the $C^\tau\sme C\rtimes_{\tau,r}H$-imprimitivity
bimodule $X$.

To form the crossed-product bimodule, we need conditions implying that
the action $\sigma: K\to\Aut C$ induces appropriately continuous
actions $\bar\sigma, \sigma$ and $\sigma\rtimes\id$ on $C^\tau, X$ and
$C\rtimes_{r,\tau}H$ respectively.  We start by observing that
$\sigma\rtimes\id:K\to\Aut (C\rtimes_{\tau,r}H)$ is always continuous.
To see this, fix $\epsilon>0$ and $f\in C_c(H,C)$.  Note that
$L:=\set{ f(s):s\in\supp f} $ is a compact subset of $C$.  Since
$\sigma:K\to\Aut C$ is continuous there exists a neighborhood $U$ of
$e_{K}$ such that $t\in U$ implies that
$\|\sigma_t(c)-c\|<\epsilon/\mu_H(\supp f)$ for all $c\in L$. Thus
$t\in U$ implies
\[
\|(\sigma\rtimes\id)_t(f)-f\|\leq\int_H\|\sigma_t\bigl(f(s)
\bigr)-f(s)\|\, ds <\epsilon.
\]

\begin{lemma}\label{lem-cts-action}
  Suppose that the map $t\mapsto \brip C\rtimes_{\tau,r}H<\sigma_t(x),
  x>$ is continuous for each fixed $x\in X_0=C_0$.  Then
  $(\bar\sigma,\sigma,\sigma\rtimes\id)$ gives a continuous action of
  $K$ on $\metrip{C^\tau}{X}{{C\rtimes}_{\tau,r}H}$.
\end{lemma}

\begin{proof}
  If $x\in X_0$ then
\begin{align*}
  \|\sigma_t(x)-x\|^2
  &=\|\brip C\rtimes_{\tau,r}H<\sigma_t(x)-x,\sigma_t(x)-x>\|\\
  &=\|(\sigma\rtimes\id)_t( \brip C\rtimes_{\tau,r}H< x, x>)
  -\brip\cctrh<\sigma_t(x), x>- \\
  &\qquad\qquad\qquad\qquad\qquad\qquad
  \brip\cctrh<x,\sigma_t(x)>+\brip\cctrh<x, x>\|\\
  &\leq\|(\sigma\rtimes\id)_t (\rip C\rtimes_{\tau,r}H< x, x>)-
  \rip\cctrh<x, x>\|+{} \\
  &\qquad\qquad \qquad\qquad
  \|\rip \cctrh<x, x>-\brip\cctrh<\sigma_t(x),x>\|+{}\\
  &\qquad\qquad\qquad\qquad\qquad\qquad\|\rip\cctrh<x,
  x>-\brip\cctrh<x,\sigma_t(x)>\|,
\end{align*}
so that $\|\sigma_t(x)-x\|\to 0$ as $t\to e_K$ using the assumption
and the continuity of $\sigma\rtimes\id$. Thus $\sigma:K\to \Aut X$ is
continuous. Since
\[
\|\bar\sigma_t(\lip {C^\tau}< x, v>)-\lip{C^\tau}< x, v>\|
\leq\|x\|\|\sigma_t(v)-v\|+\|v\|\|\sigma_t(x)-x\|
\]
the continuity of $\bar\sigma:K\to\Aut C^\tau$ follows from the
continuity of $\sigma$.
\end{proof}

Of course, by symmetry, the action $(\tau\rtimes\id,\tau,\bar\tau)$ of
$H$ on $\metrip{C\rtimes_{\sigma,r}K}{Y}{C^\sigma}$ is continuous
provided $s\mapsto\blip {C\rtimes_{\sigma,r}K}<\tau_s(y), y>$ is
continuous for each fixed $y\in Y_0=C_0$.

\begin{prop}\label{prop-tensor-ib} Suppose that the 
  action $(\bar\sigma,\sigma,\sigma\rtimes\id)$ of $K$ on
  $\metrip{C^\tau}{X}{C\rtimes_{\tau,r}H}$ is continuous and that the
  action $(\tau\rtimes\id,\tau,\bar\tau)$ of $H$ on
  $\metrip{C\rtimes_{\sigma,r}K}{Y}{C^\sigma}$ is continuous. Let $B:=
  C\rtimes_{\sigma\rtimes\tau, r}(H\times K)$. Then
  $\big(X\rtimes_{\sigma,r}K\big)\otimes_B\big(
  Y\rtimes_{\tau,r}H\big)$ is a $C^\tau\rtimes_{\bar\sigma,r}K\sme
  C^\sigma\rtimes_{\bar\tau,r}H$-imprimitivity bimodule.
\end{prop}

\begin{proof} 
  Note that whenever the actions
  $(\bar\sigma,\sigma,\sigma\rtimes\id)$ of $K$ on
  $\metrip{C^\tau}{X}{{C\rtimes}_{\tau,r}H}$ and
  $(\tau\rtimes\id,\tau,\bar\tau)$ of $H$ on
  $\metrip{C\rtimes_{\sigma,r}K}{Y}{C^\sigma}$ are continuous it makes
  sense to form the Combes $C^{\tau}\rtimes_{\bar\sigma,r}K\sme
  (C\rtimes_{\tau,r}H)\rtimes_{\sigma\rtimes\id,r}K$ and
  $(C\rtimes_{\sigma,r}K)\rtimes_{\tau\rtimes\id,r}H\sme
  C^{\sigma}\rtimes_{\bar\tau,r}H$-imprimitivity bimodules
  $X\rtimes_{\sigma,r}K$ and $Y\rtimes_{\tau,r}H$
  \cite[Remark on page 300]{com}.  We identify
  $(C\rtimes_{\tau,r}H)\rtimes_{\sigma\rtimes\id,r}K$ and
  $(C\rtimes_{\sigma,r}K)\rtimes_{\tau\rtimes\id,r}H$ with
  $B:=C\rtimes_{\sigma\rtimes\tau, r}(H\times K)$ since they are all
  naturally isomorphic. Now the internal tensor product over $B$,
  (see, for example, \cite[Proposition 3.16]{tfb}), is a
  $C^\tau\rtimes_{\bar\sigma,r}K\sme
  C^\sigma\rtimes_{\bar\tau,r}H$-imprimitivity bimodule.
\end{proof}

%%%%%%%%%%%%%%%%%%%%%%%%%%%%%%%%%%%%%%%%%%%%%%%%%%%%%%%%%%%%%%%%%%%%%%%%%%%%%%%%%%%%%%%%%%%%%%%%%%%

\section{The Concrete Imprimitivity Bimodule}
\label{concrete-result}

To get a concrete version of the imprimitivity bimodule obtained in
Proposition~\ref{prop-tensor-ib} we will use the tools developed in
\cite{hrw:xx}, where we looked at a notion of proper actions on \ib s
which generalizes Rieffel's in \cite{rie:pm88}.  (Although there are
other notions of proper actions --- for example
\cite{rie:99,exe:jfa00,mey:k00,mey:xx} --- we are closest to
\cite{rie:pm88} in spirit.)

\begin{definition}[{\cite[Definition~3.1]{hrw:xx}}]
\label{defn-pr}
If $(X,G,\gamma)$ is a Morita equivalence between two dynamical
systems $(A,G,\alpha)$ and $(B,G,\beta)$, then the action $\gamma$ of
$G$ on ${}_A X_B$ is \emph{proper} if there are an invariant subspace
$X_0$ of $X$ and invariant $*$-subalgebras $A_0$ of $A$ and $B_0$ of
$B$, such that $_{A_0}(X_0)_{B_0}$ is a pre-imprimitivity bimodule
with completion $_AX_B$, and such that
\begin{enumerate}
\item[(1)] for every $x,y\in X_0$, both $s\mapsto
  \Delta(s)^{-1/2} \blip A<x,\gamma_s(y)>$ and $s\mapsto \blip A < x ,
  \gamma_s(y)>$ are in $L^1(G,A)$;
\item[(2)] for every $b\in B_0$ and $x\in X_0$, both
  $s\mapsto \gamma_s(x)\cdot b$ and $s\mapsto
  \Delta(s)^{-1/2}\gamma_s(x)\cdot b$ are in $L^1(G,X)$;
\item[(3)]for every $x,y\in X_0$, there is a multiplier $\rip
  B^\beta<x,y>$ in $M(B_0)^\beta$ such that $z\cdot \rip
  B^\beta<x,y>\in X_0$ for all $z\in X_0$, and
\begin{equation}\label{pr-eq}
\int_G b\beta_s(\rip B<x,y>)\, ds
=b\rip B^\beta<x,y>\quad\text{for all $b\in B_0$}.
\end{equation}
\end{enumerate}  
\end{definition}
If $\gamma$ is also saturated, as in \cite[Definition~3.15]{hrw:xx},
then $(X_0)_{B^\beta}$ completes to give an imprimitivity bimodule
implementing a Morita equivalence between $A\rtimes_{\alpha,r}G$ and a
generalized fixed-point algebra $B^\beta$ of $B$ \cite[Theorem
3.16]{hrw:xx}.  Furthermore, in \cite[Theorem~4.1]{hrw:xx} we showed
that the action $\beta$ on $B$ is proper with respect to $B_{1} := \rip B<X_{0},X_{0}>$, so $B_1$ completes to a
$B\rtimes_{\beta,r}G\sme B^\beta$-\ib. Finally, a linking algebra argument from
\cite{hrw:xx} implies that there is a natural \ib-isomorphism of
\begin{equation*}
  (X\rtimes_{\alpha,r}G)\tensor_{B\rtimes_{\beta,r}G}\overline{B_{1}}
  \quad \text{onto}\quad \overline{X_{0}}.
\end{equation*}
 Note that if
${}_{A_0}(X_0)_{B_0}={}_{B_0}(B_0)_{B_0}$ then
Definition~\ref{defn-pr} reduces to that of Rieffel's and
\cite[Theorem 3.16]{hrw:xx} reduces to \cite[Corollary
1.7]{rie:pm88}

In our situation, we want the action $\sigma$ of $K$ on
$\metrip{C^\tau}{X}{C\rtimes_{\tau,r}H}$ to be proper with respect to
$\metrip{(D_0,\bar\sigma)}{(X_0,\sigma)}{(E_0,\sigma\rtimes\id)}$.  We
will then obtain the concrete version of the symmetric imprimitivity
theorem, as well as the desired isomorphism onto the tensor-product
version, along the following lines: if the action on $X$ is proper and
saturated, then first, \cite[Theorem~3.16]{hrw:xx} implies that
$\overline{X_{0}}$ is a $C^{\tau}\rtimes_{\bar\sigma,r}K\sme
(C\rtimes_{\tau,r}H)^{\sigmartimesid}$-\ib.  On the other hand,
\cite[Theorem~4.1]{hrw:xx} implies $\sigmartimesid$ is saturated and
proper with respect to $E_{0}\subset C\rtimes_{\tau,r}H$ so that
$\overline{E_{0}}$ is a
$(C\rtimes_{\tau,r}H)\rtimes_{\sigmartimesid,r}K \sme
(C\rtimes_{\tau,r}H)^{\sigmartimesid}$-\ib.  Let $Y$ be the
$C\rtimes_{\sigma, r}K\sme C^\sigma$-imprimitivity bimodule based on
$Y_0:=C_0$ coming from the proper and saturated action $\sigma$ of $K$
on $C$.  Then the Combes bimodule $Y\rtimes_{\tau,r}H$ is a
$(C\rtimes_{\sigma,r}K)\rtimes_{\sigma\times \id}H\sme
C^{\sigma}\rtimes _{\bar\tau,r}H$-\ib.  We can identify
$(C\rtimes_{\sigma,r}K)\rtimes_{\sigma\times \id}H $ with
$C\rtimes_{\tau\times \sigma}(K\times H)$, and, with a bit more work,
we will show that $(C\rtimes_{\tau,r}H)^{\sigmartimesid}$ and
$C^{\sigma}\rtimes _{\bar\tau,r}H$ are isomorphic.  With these
identifications, we will show that $\overline{E_{0}}$ is isomorphic to
$Y\rtimes_{\tau,r}H$, and then \cite[Theorem~4.1]{hrw:xx} gives an
isomorphism
\begin{equation*}
  (X\rtimes_{\sigma,r}K)\tensor_{C\rtimes_{\sigma\times \tau,r}(K\times
  H)} (Y\rtimes_{\tau,r}H) \cong \overline{X_0}
\end{equation*}
as $C^{\tau}\rtimes_{\bar\sigma,r}K\sme
C^{\sigma}\rtimes_{\bar\tau}H$-\ib s.  This gives us both the desired
concrete Morita equivalence and the isomorphism onto the
tensor-product version.

Thus, our first step is to find conditions which ensure that the three
items of Definition~\ref{defn-pr} hold in our situation so that the
action $\sigma$ of $K$ on $\metrip{C^\tau}{X}{C\rtimes_{\tau,r}H}$ is
proper with respect to
$\metrip{(D_0,\bar\sigma)}{(X_0,\sigma)}{(E_0,\sigma\rtimes\id)}$.

For Definition~\ref{defn-pr}(1) we need to assume that for every
$x,v\in X_0$ the functions
\begin{equation}\label{eq-left-integrable}
t\mapsto\blip{C^\tau}<x,\sigma_t(v)>
\quad\text{and}\quad t\mapsto\Delta_K(t)^{-1/2}
\blip{C^\tau}<x,\sigma_t(v)>
\end{equation}
are in $L^1(K, C^\tau)$.  Note that if $w\in X_0$ then
$\sigma_t(x)\cdot\rip \cctrh<v, w>=\blip {C^\tau}<\sigma_t(x),v>\cdot
w$, so that the integrability of the functions in
\eqref{eq-left-integrable} implies the integrability of $t\mapsto
\sigma_t(x)\cdot e$ and its product with $\Delta_K(t)^{-1/2}$ for all
$x\in X_0$ and $e\in E_0$.  So, Definition~\ref{defn-pr}(1) and~(2)
hold provided the functions in \eqref{eq-left-integrable} are
integrable. For Definition~\ref{defn-pr}(3), we need to assume that
the function $s\mapsto\Delta_H(s)^{-1/2}\brip C^{\sigma}<x,\tau_s(v)>$
is in $L^1(H,C^\sigma)$ whenever $x,v\in X_0$ and, using the
$C^{\sigma}$-valued inner product for $Y$, define
\begin{equation}\label{eq-ripmult}
\brip(C\rtimes_{\tau,r}H)^{\sigma\rtimes\id}<x,
v>(s):=\Delta_H(s)^{-1/2}\brip C^\sigma<x,\tau_s(v)>;
\end{equation}
now we still need to find conditions which ensure that
\begin{enumerate}
\item[(1)] $\rip(C\rtimes_{\tau,r}H)^{\sigma\rtimes\id}<x, v>$ is a
  multiplier of $C\rtimes_{\tau,r}H$;
\item[(2)] $\rip(C\rtimes_{\tau,r}H)^{\sigma\rtimes\id}<x, v>$
  multiplies $E_0$ and is invariant under $\sigma\rtimes\id$;
\item[(3)] equation \eqref{pr-eq} is satisfied; and
\item[(4)] $w\cdot\rip(C\rtimes_{\tau,r}H)^{\sigma\rtimes\id}<x, v>$
  is back in $X_0$.
\end{enumerate}

Item~(1) is an immediate consequence of Proposition~\ref{prop-mult}
below which allows us to view elements of $L^{1}(H,C^{\sigma})$ as
multipliers of $C\rtimes_{\sigma}H$ via convolution.  Although
essentially straightforward and presumably known, 
its proof requires some intricacies of
vector-valued integration which are
certainly far from the subject at hand.  We provide a detailed proof in
Appendix~\ref{appendix1}.

\begin{prop}\label{prop-mult}
  Let $(A,G,\alpha)$ be a dynamical system.  Suppose that $B$ is a
  subalgebra of $M(A)$ such that $(B,G,\bar\alpha)$ is a dynamical
  system.  If $g\in L^{1}(G,B)$, then there is a unique multiplier
  $T_{g}$ in $M(A\rtimes_{\alpha}G)$ such that for all $f\in
  L^{1}(G,A)$ both $T_{g}f$ and $fT_{g}$ are in $L^{1}(G,A)$ (viewed
  as a subalgebra of $A\rtimes_{\alpha}G$), and for almost all $s$,
  \begin{align}
    T_{g}f(s) &= \int_{G} g(r) \alpha_{r}\bigl(f(r^{-1}s)\bigr)
    \,dr\quad \text{and} \label{eq:5}\\
    fT_{g}(s) &= \int_{G} f(r) \bar\alpha_{r}\bigl(g(r^{-1}s)\bigr)
    \,dr. \label{eq-prop-mult}
  \end{align}
\end{prop}

Note that if $B=A$, then Proposition~\ref{prop-mult} reduces to two
things. First, the familiar formula for convolutions of functions in
$C_c(G,A)$ extends to functions in $L^1(G,A)$, and second, that
convolution has the usual $*$-algebraic properties.

For (2), we seek conditions which ensure that
$\brip(C\rtimes_{\tau,r}H)^{\sigma\rtimes\id}<x, v>$ multiplies $E_0$.
Using \eqref{eq-prop-mult} we compute:
\begin{align*}
  \rip C\rtimes_{\tau,r}H< u,
  w>&\rip(C\rtimes_{\tau,r}H)^{\sigma\rtimes\id} <x, v>(s)
  \\
  &= \int_H \rip\cctrh<u, w>(r)\bar\tau_r\big(
  \rip(C\rtimes_{\tau,r}H)^{\sigma\rtimes\id}<x, v>(r^{-1}s)\big)\, dr\\
  &=\Delta_H(s)^{-1/2}u^*\int_H\tau_r(w)\bar\tau_r\big(
  \brip C^\sigma< x, \tau_{r^{-1}s}(v)> \big)\, dr\\
  &=\Delta_H(s)^{-1/2}u^*\tau_s\Bigl( \int_H \tau_{s^{-1}r}(w)
  \brip C^\sigma< \tau_{s^{-1}r}(x), v>\, dr \Bigr)\\
  &=\Delta_H(s)^{-1/2}u^*\tau_s\Bigl( \int_H \tau_{r}(w) \brip
  C^\sigma<\tau_{r}(x), v>\, dr \Bigr)
\end{align*}
by the change of variable $s^{-1}r\mapsto r$.  Since
\begin{align*}
  \int_H \tau_{r}(w)\brip C^\sigma< \tau_{r}(x), v>\, dr &=\int_H
  \blip {C\rtimes_{\sigma,r}K}< \tau_{r}(w),\tau_{r}(x)>
  \cdot v\, dr\\
  &=\int_H\int_K \blip {C\rtimes_{\sigma,r}K}<
  \tau_{r}(w),\tau_{r}(x)>(t)
  \sigma_t(v)\Delta_K(t)^{1/2}\, dt\, dr\\
  &=\int_H\int_K\tau_r(w)\sigma_t(\tau_r(x^*)v)\, dt\, dr.
\end{align*}
Thus $\rip C\rtimes_{\tau,r}H< u, w>\rip
(C\rtimes_{\tau,r}H)^{\sigma\rtimes\id} <x, v>$ is back in $E_0$
provided
\begin{equation}\label{eq-back-x0}
\int_H\int_K\tau_r(w)\sigma_t(\tau_r(x^*)v)\, dt\, dr \in X_0
\end{equation} 
whenever $u,v,w,x\in X_0$.

For (3), we need to show that
\begin{equation*}
\int_K e(\sigma\rtimes\id)_t(\rip C\rtimes_{\tau,r}H<
x, v>\, )dt=e\rip (C\rtimes_{\tau,r}H)^{\sigma\rtimes\id}<x,
v>.
\end{equation*}
To do this we need the following lemma. For the sake of clarity, we
have decorated our integrals with the space in which the integral
takes values; thus we write $\int_H^C f(s)\,
ds$ for the $C$-valued integral of $f\in L^1(H,C)$.  Again, the proof of the lemma
requires some gymnastics with vector-valued integration, so we
relegate it to Appendix~\ref{appendix2}.

\begin{lemma}\label{lem-variable}
  Assume that for all $u,v,w,x\in X_{0}$, the function
  \begin{equation*}
    (r,s,t)\mapsto u\tau_{r}(v^{*})\sigma_{t}\bigl( \tau_{r}(w^{*})
    \tau_{s}(x) \bigr)\Delta_{H}(s)^{-\frac12}
  \end{equation*}
  is integrable as a function from $H\times H\times K$ to $C$.  Then
  for all $e_{i}\in E_{0}$ the function $t\mapsto
  e_{1}(\sigma\rtimes\id)_{t}(e_{2})$ is integrable as a function from
  $K$ to $C\rtimes_{\tau,r}H$.  Further, the integral
  \begin{equation}\label{eq:6}
    \int_{K}^{C\rtimes_{\tau,r}H}e_{1}(\sigma\rtimes\id)_{t}(e_{2})\, dt
  \end{equation}
  takes values in $L^1(H,C)$ viewed as a subalgebra of
  $C\rtimes_{r,\tau}H$, and a representative for \eqref{eq:6} is given
  by
\begin{equation}  \label{eq:8}
s\mapsto \int_{K}^{C}e_{1}(\sigma\rtimes\id)_{t}(e_{2})(s)\, dt.
\end{equation}
\end{lemma}

Using Lemma~\ref{lem-variable} we obtain
\begin{align*}
  \int_K e(\sigma\rtimes\id)_t&(\rip C\rtimes_{\tau,r}H <x, v>)\,
  dt(s) =\int_K e(\sigma\rtimes\id)_t
  (\brip C\rtimes_{\tau,r}H< x, v>)(s)\, dt \\
  &=\int_K\int_H e(r)\tau_r\bigl( (\sigma\rtimes\id)_t(\brip
  C\rtimes_{\tau,r}H < x, v>(r^{-1}s)
  \bigr)\, dr\, dt\\
  &=\int_K\int_H e(r) \tau_r\sigma_t(x^*\tau_{r^{-1}s}(v))
  \Delta_H(r^{-1}s)^{-1/2}\, dr\, dt  \\
  \intertext{which, since our assumptions guarantee that the integrand
    is an integrable function of $(r,t)$ provided $e\in E_{0}$ (see
    Remark~\ref{rem-three-fcns}), is, by Fubini's Theorem, equal to}
  &=\int_H\int_K e(r)\sigma_t\bigl(
  \tau_r(x)^*\tau_s(v) \bigr)\Delta_H(r^{-1}s)^{-1/2}\, dt\, dr \\
  &=\int_H e(r) \brip C^\sigma<\tau_r(x),\tau_s(v)>
  \Delta_H(r^{-1}s)^{-1/2}\, dr  \\
  &=\int_H e(r)\bar\tau_r\bigl( \rip
  C^\sigma<x,\tau_{r^{-1}s}(v)>\bigr )
  \Delta_H(r^{-1}s)^{-1/2}\, dr  \\
  &=\int_H e(r)\bar\tau_r\bigl( \rip
  (C\rtimes_{\tau,r}H)^{\sigma\rtimes\id}<x, v>(r^{-1}s) \big)\,
  dr \\
  &=e\rip (C\rtimes_{\tau,r}H)^{\sigma\rtimes\id}<x, v>(s)
  \end{align*}
  as required.
  
  To establish (4), we claim that
\begin{equation}\label{eq:10}
  w\cdot{}\rip (C\rtimes_{\tau,r}H)^{\sigma\rtimes\id}<x,
  v>=\int_K
  w\cdot(\sigma\rtimes\id)_t(\rip C\rtimes_{\tau,r}H<x,
  v>)\, dt.
\end{equation}
To see this, first note that
\begin{align*}
  \|w\cdot (\sigmartimesid)_{t}\bigl(\rip\cctrh<x,v>\bigr) \| &= \|
  w\cdot \brip\cctrh <\sigma_{t}(x),\sigma_{t}(v)>\| \\
  &\le \|\blip
  C^{\tau}<w, \sigma_{t}(x)>\|\|\sigma_{t}(v)\| \\
  &\le \| \blip C^{\tau}<w, \sigma_{t}(x)> \|\|v\|,
\end{align*}
and the latter is integrable by Definition~\ref{defn-pr}(1).  It is
not hard to see that
\begin{equation*}
  \brip \cctrh <z,\int_{K} w\cdot (\sigmartimesid)_{t}\bigl({ \rip
  \cctrh<x,v>}\bigr) \,dt >=\brip \cctrh <z, w\cdot
  {\rip(\cctrh)^{\sigmartimesid} <x,v>}>,
\end{equation*}
and \eqref{eq:10} follows. We use \eqref{dual-ra} to write the
left-hand side of \eqref{eq:10} as
\begin{align*}
  \int_K\int_H \tau_s^{-1}\Bigl(w\bigl({}&{}(\sigma\rtimes\id)_t (\rip
  C\rtimes_{\tau,r}H<x, v>) (s)\bigr)\Bigr)
  \Delta_H(s)^{-1/2}\, ds\, dt \\
  &= \int_K\int_H\tau_s^{-1}\Bigl( w\sigma_t\bigl( x^*\tau_s(v) \bigr)
  \Bigr)
  \Delta_H(s)^{-1}\, ds\, dt \\
  &= \int_K\int_H\tau_s^{-1}(w)
  \sigma_t\bigl(\tau_{s^{-1}}(x^*)v \big)\Delta_H(s)^{-1}\,ds\, dt  \\
  &= \int_K\int_H \tau_s(w)\sigma_t(\tau_s(x^*)v)\, ds\, dt,
\end{align*}
and we already assumed in Equation~\eqref{eq-back-x0} above that this
double integral is back in $X_0$.
 
We can restate our conclusions above as

\begin{prop}\label{prop-symmetric-proper}
  Suppose that the action $(\bar\sigma,\sigma,\sigma\rtimes\id)$ of
  $K$ on $\metrip{C^\tau}{X}{C\rtimes_{\tau,r}H}$ is continuous. If
  for all $u,v,w,x\in X_0$
\begin{enumerate}
\item[(1)] the function $t\mapsto\blip {C^\tau}<x,\sigma_t(v)>$ and
  its product with $\Delta_K(t)^{-1/2}$ are in $L^1(K, C^\tau)$;
\item[(2)] the function $s\mapsto\Delta_H(s)^{-1/2} \brip
  C^\sigma<x,\tau_s(v)>$ is in $L^1(H,C^\sigma)$;
\item[(3)] the integral $\int_H\int_K\tau_r(w)\sigma_t(\tau_r(x^*)v)\,
  dt\, dr$ is in $X_0$; and
\item[(4)] the function $(r,s,t)\mapsto u\tau_r(v^*)\sigma_t\big(
  \tau_r(w^*)\tau_s(x) \big)\Delta_H(s)^{-1/2}$ is integrable,
\end{enumerate}
then $\sigma$ is a proper action of $K$ on
$\metrip{C^\tau}{X}{C\rtimes_{\tau,r}H}$ with respect to
$\metrip{(D_0,\bar\sigma)}{(X_0,\sigma)}{(E_0,\sigma\rtimes\id)}$.
\end{prop}

In the situation of Proposition~\ref{prop-symmetric-proper} above, we
want the action of $K$ to be saturated with respect to $X_0$, so that
$X_0$ completes to a $C^\tau\rtimes_{\bar\sigma,r}K\sme
(C\rtimes_{\tau,r}H)^{\sigma\rtimes\id}$-imprimitivity bimodule by
\cite[Theorem~3.16]{hrw:xx}.  Since $\sigma$ is a proper action of $K$
on $\metrip{C^\tau}{X}{C\rtimes_{\tau,r}H}$ with respect to
$\metrip{(D_0,\bar\sigma)}{(X_0,\sigma)}{(E_0,\sigma\rtimes\id)}$,
\cite[Theorem~4.1]{hrw:xx} says that $\sigma\rtimes\id$ is a proper
action on $C\rtimes_{\tau,r}H$ with respect to $\rip\cctrh <X_0,
X_0>$. We have set things up so that $\brip C\rtimes_{\tau,r}H< X_0,
X_0>=E_0$, and thus $E_0$ completes to an $I$--$J$-imprimitivity
bimodule, where $I$ is an ideal in
$(C\rtimes_{\tau,r}H)\rtimes_{\sigma\rtimes\id}K$ and $J$ is a
generalized fixed point algebra of $C\rtimes_{\tau,r}H$, by
\cite[Theorem 1.5]{rie:pm88}.  But \cite[Theorem~4.1]{hrw:xx} implies
that $J$ is an ideal in $(C\rtimes_{\tau,r}H)^{\sigma\rtimes\id}$, and
if $\sigma\rtimes\id$ is saturated with respect to $E_0$ then $\sigma$
is saturated with respect to $X_0$, and
$J=(C\rtimes_{\tau,r}H)^{\sigma\rtimes\id}$.
In applications we expect that it will be easy to check that
$\sigma\rtimes\id$ is saturated with respect to $E_0$, and then we
have

\begin{prop}
  Suppose that the action $(\bar\sigma,\sigma,\sigma\rtimes\id)$ of
  $K$ on $\metrip{C^\tau}{X}{C\rtimes_{\tau,r}H}$ is continuous and
  proper with respect to
  $\metrip{(D_0,\bar\sigma)}{(X_0,\sigma)}{(E_0,\sigma\rtimes\id)}$,
  and that the proper action $\sigma\rtimes\id$ of $K$ on
  $C\rtimes_{\tau,r}H$ is saturated with respect to $E_0$. Then $E_0$
  completes to a
  $(C\rtimes_{\tau,r}H)\rtimes_{\sigma\rtimes\id,r}K\sme
  (C\rtimes_{\tau,r}H)^{\sigma\rtimes\id}$-imprimitivity bimodule.
\end{prop}

In the next proposition we identify the module $\overline{E_0}$ and
the fixed point algebra $(C\rtimes_{\tau,r}H)^{\sigma\rtimes\id}$;
again we add substantial hypotheses.  The set $E_0$ is always a subset
of $L^1(H, C_0)$.  Even though the completion $Y$ of $C_{0}$ is a
$C\rtimes_{\sigma,r}K\sme C^{\sigma}$-\ib, this does not in general
imply that $E_0$ is contained in $L^1(H, Y)$ (the norm of $Y$ is not
related to the norm of $C$).

\begin{prop}\label{prop-saturation}
  Suppose that the action $(\tau\rtimes\id,\tau,\bar\tau)$ of $H$ on
  $\metrip{C\rtimes_{\sigma,r}K}{Y}{C^\sigma}$ is continuous so that
  $Y\rtimes_{\tau,r}H$ is a $(C\rtimes_{\sigma,r}K)
  \rtimes_{\tau\rtimes\id,r}H\sme
  C^\sigma\rtimes_{\bar\tau,r}H$-imprimitivity bimodule. Also suppose
  that the action $\sigma\rtimes\id$ of $K$ on $C\rtimes_{\tau,r}H$ is
  proper and saturated with respect to $E_0$.  If $E_0\subset L^1(H,
  Y)$ and if for all $u,v,w,x\in X_0$ the function
%  $(r,s,t)\mapsto u\tau_r(v^*)\sigma_t\big( \tau_r(w^*)\tau_s(x)
%  \big)\Delta_H(s)^{-1/2}$ 
  given in Proposition~\ref{prop-symmetric-proper}(4) is integrable,
  then $(C\rtimes_{\tau,r}H)^{\sigma\rtimes\id}\cong
  C^\sigma\rtimes_{\bar\tau,r}H$ and
\[
{}_{(C\rtimes_{\tau,r}H) \rtimes_{\sigma\rtimes\id}K}
\overline{E_0}_{(C\rtimes_{\tau,r}H)^{\sigma\rtimes\id}} \cong
{}_{(C\rtimes_{\sigma,r}K) \rtimes_{\tau\rtimes \id,r}H}
(Y\rtimes_{\tau,r}H)_{C^\sigma\rtimes_{\bar\tau,r}H}
\]
as imprimitivity bimodules.
 \end{prop}

\begin{proof}
  As observed in \cite[page 14]{ekqr:bams00}, if $h\in C_c(H,Y)$ then
  the Cauchy-Schwartz inequality gives $\|h\|\leq\int_H \|h(s)\|\,
  ds=\|h\|_1$. It follows that $L^1(H,Y)$ is dense in
  $Y\rtimes_{\tau,r}H$ and the actions and inner products on
  $C_{c}(H,Y)\subset Y\rtimes_{\tau,r}H$ given, for example, in
  equations (4.6)--(4.9) in \cite{hrw:xx},
%  \eqref{combesla}--\eqref{combesrip} 
  extend to integrable functions.  That the formulas themselves extend
  to integrable functions can be seen by including
  $Y\rtimes_{\tau,r}H$ in the linking algebra $L(Y\rtimes_{\tau,r}H)$
  and doing the computations in the relevant bit of the $C^*$-algebra.
  
  By assumption $E_0$ is contained in $L^1(H,Y)$; we let $\iota:E_0\to
  Y\rtimes_{\tau,r}H$ be the inclusion map.  Let
  $\phi:(C\rtimes_{\tau,r}H)\rtimes_{\sigma\rtimes\id,r}K\to
  (C\rtimes_{\sigma,r}K)\rtimes_{\tau\rtimes\id,r}H$ be the
  isomorphism such that $\phi(f)(s)(t)=f(t)(s)$ for $f\in
  L^1(K,L^1(H,C))$. We will show that $\rip
  (C\rtimes_{\tau,r}H)^{\sigma\rtimes\id}< e, f> \mapsto \rip
  C^\sigma\rtimes_{\bar\tau,r}H<e, f>$ extends to an isomorphism
  $\psi$ of $(C\rtimes_{r,\tau}H)^{\sigma\rtimes\id}$ onto
  $C^\sigma\rtimes_{\bar\tau}H$ and that $(\phi,\iota,\psi)$ extends
  to an imprimitivity bimodule isomorphism of $\overline{E_0}$ onto
  $Y\rtimes_{\tau,r}H$.
  
  If $e,f\in E_0$ then
\begin{align}
  \lip {(C\rtimes_{\tau,r}H)\rtimes_{\sigma\rtimes\id,r}K}<e, f>(t)(s)
  &=\Delta_K(t)^{-1/2}\blip
  {C\rtimes_{\tau,r}H}< e,(\sigma\rtimes\id)_t(f)> (s)\notag\\
  &=\Delta_K(t)^{-1/2} e*(\sigma\rtimes\id)_t(f)^*(s)\notag\\
  &=\Delta_K(t)^{-1/2}\int_H
  e(r)\tau_r((\sigma\rtimes\id)_t(f^*)(r^{-1}s))\, dr\notag\\
  &=\Delta_K(t)^{-1/2}\int_H
  e(r)\tau_s\sigma_t(f(s^{-1}r)^*)\Delta_H(s^{-1}r)\, dr\notag\\
  &=\int_H \Delta_K(t)^{-1/2}
  e(r)\sigma_t(\tau_s(f(s^{-1}r))^*)\Delta_H(s^{-1}r)\, dr\notag\\
  &=\int_H \blip {C\rtimes_{\sigma,r}K}
  < e(r),\tau_s(f(s^{-1}r))>(t)\Delta_H(s^{-1}r)\, dr\notag\\
  &=\lip{(C\rtimes_{\sigma,r}K)\rtimes_{\tau\rtimes\id,r}H}< e,
  f>(s)(t)\notag
\end{align}
so that $\phi(\lip{(C\rtimes_{\tau,r}H)\rtimes_{\sigma\rtimes\id,r}K}<
e, f>)=\blip {(C\rtimes_{\sigma,r}K)\rtimes_{\tau\rtimes\id,r}H}<
\iota(e), \iota(f)>$.  Next, we show that
\[
\iota\bigl( \lip{(C\rtimes_{\tau,r}H)\rtimes_{\sigma\rtimes\id,r}K}<e,
f> \cdot g\big) = \phi(\lip
{(C\rtimes_{\tau,r}H)\rtimes_{\sigma\rtimes\id,r}K} <e, f>)\cdot
\iota(g)
\]
whenever $e,f,g\in E_0$. First,
\begin{align}
  \lip {(C\rtimes_{\tau,r}H)\rtimes_{\sigma\rtimes\id,r}K}<e, f>\cdot
  g (s) &=
  e\cdot \rip (C\rtimes_{\tau,r}H)^{\sigma\rtimes\id}<f, g>(s) \notag\\
  &=\int_K e*((\sigma\rtimes\id)_t(f^**g))\,dt(s)\notag\\
  \intertext{which, by Lemma~\ref{lem-variable}, is}
  &=\int_K e*((\sigma\rtimes\id)_t(f^**g))(s)\,dt \notag\\
  &=\int_K\int_H e(r)\tau_r\big( (\sigma\rtimes\id)_t(f^**g)(r^{-1}s)  \big)\, dr\, dt\notag\\
  &=\int_K\int_H\int_H e(r)\sigma_t\tau_{ru^{-1}}\big(
  f(u)^*g(ur^{-1}s) \big)\, du\, dr\, dt\label{eq-lh-map}.
\end{align}
On the other hand,
\begin{align*}
  \phi(\lip {(C\rtimes_{\tau,r}H)\rtimes_{\sigma\rtimes\id,r}K} <e,
  f>){}&{}\cdot \iota(g)= \lip
  {(C\rtimes_{\sigma,r}K)\rtimes_{\tau\rtimes\id,r}H}<e,
  f>\cdot g (s) \\
  &= e\cdot\brip C^\sigma\rtimes_{\bar\tau,r}H< f, g>(s)\\
  \intertext{which, using the Combes action, is} &= \int_H
  e(r)\cdot\bar\tau_r\bigl(
  \rip C^\sigma\rtimes_{\bar\tau,r}H<f, g>(r^{-1}s) \bigr)\, dr\\
  &= \int_H e(r)\bar\tau_r\Bigl( \int_H\bar\tau_u^{-1} \bigl( \brip
  C^\sigma< f(u), g(ur^{-1}s)> \bigr)
  \Bigr)\, du\,dr\\
%  &=\int_H\int_H e(r)
%  \brip C^\sigma< \tau_{ru^{-1}}(f(u)),\tau_{ru^{-1}}(g(ur^{-1}s)) >\,
%  du
%\, dr\\
  &=\int_H\int_H\int_K e(r)\sigma_t\tau_{ru^{-1}}\big(
  f(u)^*g(ur^{-1}s)\big)\, dt\,du\,dr
\end{align*}
which is the same as~\eqref{eq-lh-map} by two applications of
Fubini's Theorem.

So far we have shown that
\begin{equation*}
(\phi,\iota):
\metrip{(C\rtimes_{\tau,r}H)\rtimes_{\sigma\rtimes\id,r}K}{E_0}{}\to
\metrip{(C\rtimes_{\sigma,r}K)\rtimes_{\tau\rtimes\id,r}H}
{(Y\rtimes_{\tau,r}H)}{}
\end{equation*}
extends to a monomorphism of left-Hilbert modules. To produce
$\psi:(C\rtimes_{r,\tau}H)^{\sigma\rtimes\id} \to
C^\sigma\rtimes_{\bar\tau}H$ taking
\begin{equation*}
\rip (C\rtimes_{\tau,r}H)^{\sigma\rtimes\id}<e,
f>\quad\text{to}\quad \brip C^\sigma\rtimes_{\bar\tau,r}H<\iota(e),
\iota(f)>,
\end{equation*}
let $\pi$ be a faithful representation of $C$ on $\H_\pi$ and
$\bar\pi$ its extension to $M(C)$. By our assumptions,
$\bar\tau:H\to\Aut C^\sigma$ is continuous, so that the regular
representation $(\widetilde\pi,\lambda)$ of $(C,H,\tau)$ extends to a
covariant representation $(\widetilde{\bar\pi},\lambda)$ of
$(C^\sigma,H,\bar\tau)$. The representations
\begin{gather*}
  \overline{\widetilde\pi\rtimes\lambda}:
  (C\rtimes_{\tau,r}H)^{\sigma\rtimes\id}\to\B(L^2(H,\H_\pi))\quad \text{and} \\
  \widetilde{\bar\pi}
  \rtimes\lambda:C^\sigma\rtimes_{\bar\tau,r}H\to\B(L^2(H,\H_\pi))
\end{gather*}
are faithful. We will show that for $e,f\in E_0$
\begin{equation}
\overline{\widetilde\pi\rtimes\lambda}
(\rip (C\rtimes_{\tau,r}H)^{\sigma\rtimes\id}<e,
f>)=\widetilde{\bar\pi}\rtimes\lambda(\brip
C^\sigma\rtimes_{\bar\tau,r}H
<\iota(e), \iota(f)>).
\end{equation}
It then follows that
$\psi:=(\widetilde{\bar\pi}\rtimes\lambda)^{-1}\circ
\overline{\widetilde\pi\rtimes\lambda}$ is an injective homomorphism
satisfying
\begin{equation}
  \label{eq:4}
  \psi\bigl(\rip (C\rtimes_{\tau,r}H)^{\sigma\rtimes\id}<e,
f>\bigr) =\brip
C^\sigma\rtimes_{\bar\tau,r}H
<\iota(e), \iota(f)>.
\end{equation}
% and
% \begin{equation}
% \psi(\rip (C\rtimes_{\tau,r}H)^{\sigma\rtimes\id}<e, f>(s)=
% \brip C^\sigma <\iota(e), \tau_s(\iota(f))>
% \Delta_H(s)^{-1/2}.\label{eq-psi-mult}
% \end{equation}
As the first step we show that $\iota\bigl(g\cdot \rip
(C\rtimes_{\tau,r}H)^{\sigma\rtimes\id}<e, f>\bigr)=\iota(g)\cdot\brip
(C\rtimes_{\tau,r}H)^{\sigma\rtimes\id} <\iota(e), \iota(f)>$ whenever
$e,f,g\in E_0$; this will then also give us that $(\iota,\psi)$
extends to a right-Hilbert module homomorphism.  Well,
\begin{align*}
  \iota\big( g\cdot\rip (C\rtimes_{\tau,r}H)^{\sigma\rtimes\id} <e,
  f>\big) &=\iota\big( \lip
  {(C\rtimes_{\tau,r}H)\rtimes_{\sigma\rtimes\id,r}K}<g,
  e>\cdot f \big) \\
  &=\phi\big( \lip
  {(C\rtimes_{\tau,r}H)\rtimes_{\sigma\rtimes\id,r}K}<g, e>
  \big)\cdot \iota(f)\\
  &=\blip {(C\rtimes_{\sigma,r}K)\rtimes_{\tau\rtimes\id,r}H}<
  \iota(g),\iota(e)>\cdot\iota(f) \\
  &=\iota(g)\cdot\brip
  C^\sigma\rtimes_{\bar\tau,r}H<\iota(e),\iota(f)>.
\end{align*}

Note that for all $g\in E_0$ and $\xi,\eta\in L^2(H,\H_\pi)$
\begin{align*}
  \bip(\overline{\widetilde\pi\rtimes\lambda}( \rip
  (C\rtimes_{\tau,r}H)^{\sigma\rtimes\id}<e,f>)\xi |{}&{}
  \widetilde\pi\rtimes\lambda(g^*)\eta ) =\bip (
  \widetilde\pi\rtimes\lambda(g\cdot
  \rip(C\rtimes_{\tau,r}H)^{\sigma\rtimes\id}<e,f>)\xi  |\eta ) \\
  &=\bip( \widetilde\pi\rtimes\lambda(\iota(g)\cdot \brip
  C^\sigma\rtimes_{\bar\tau,r}H< \iota(e),\iota(f)>)\xi |
  \eta) \\
  &=\bip(\widetilde{\bar\pi}\rtimes\lambda( \brip
  C^\sigma\rtimes_{\bar\tau,r}H < \iota(e),\iota(f)>)\xi
  |\widetilde\pi\rtimes\lambda(g^*)\eta ),
\end{align*}
so that \eqref{eq:4} follows from the nondegeneracy of
$\widetilde\pi\rtimes\lambda$.

We now know that the triple $(\phi,\iota,\psi)$ is an isomorphism of
$\overline{E_0}$ into $Y\rtimes_{\tau,r}H$; to see that this map is
onto $Y\rtimes_{\tau,r}H$ we will show that $\iota(E_0)$ is invariant
under the right action of $C^\sigma$ and $H$, and then must be
invariant under the action of $C^\sigma\rtimes_{\bar\tau,r}H$ as well.
This suffices because $\phi$ maps
$(C\rtimes_{\tau,r}H)\rtimes_{\sigma\rtimes\id,r}K$ onto
$(C\rtimes_{\sigma,r}K)\rtimes_{\tau\rtimes\id,r}H$ and the Rieffel
correspondence then implies that
$\overline{\iota(E_0)}=Y\rtimes_{\tau,r}H$.  Therefore $\psi:
(C\rtimes_{r,\tau}H)^{\sigma\rtimes\id}\to
C^\sigma\rtimes_{\bar\tau,r}H$ is an isomorphism.

Consider $e\in E_0$, say $e(s)=\Delta_H(s)^{-1/2}x\tau_s(y^*)$, and
$\brip C^\sigma< v, w>$ where $v,w,x,y\in X_0$. Then
\begin{align*}
  (e\cdot \brip C^\sigma< v, w>)(s)=e(s)\tau_s(\rip C^\sigma<v, w>)
  &=\Delta_H(s)^{-1/2}x\tau_s(y^*)\tau_s(\langle v\,,\,
  w\rangle_{C^\sigma})\\&=\Delta_H(s)^{-1/2}x\tau_s(y^*\rip
  C^\sigma<v, w>)
\end{align*}
so that $e\cdot \brip C^\sigma< v, w>$ is back in $E_0$ because $\brip
C^\sigma< v, w>$ multiplies $X_0$.  Also,
\[
(e \cdot i_H(r))(s)=e(sr^{-1})\Delta_H(r^{-1})
=\Delta_H(sr^{-1})^{-1/2}\bigl(\Delta_H(r^{-1})x\tau_{sr^{-1}}(y^*)\bigr)
\]
so $e \cdot i_H(r)$ is back in $E_0$.
\end{proof}

We can summarize our above discussions as follows.
\begin{thm}\label{thm-symmetric}
  Suppose that $\tau:H\to\Aut C$ and $\sigma: K\to\Aut C$ are
  commuting actions which are proper and saturated with respect to the
  same dense $*$-subalgebra $C_0$.
\begin{enumerate}
\item[(1)](Continuity) If the maps $t\mapsto \brip
  C\rtimes_{\tau,r}H<\sigma_t(x), x>$ and $s\mapsto \blip
  {C\rtimes_{\sigma,r}K}<\tau_s(x), x>$ are continuous for all $x\in
  X_0=Y_0=C_0$ then the actions $(\bar\sigma,\sigma,\sigma\rtimes\id)$
  of $K$ on $\metrip{C^\tau}{X}{C\rtimes_{\tau,r}H}$ and
  $(\tau\rtimes\id,\tau,\bar\tau)$ of $H$ on
  $\metrip{C\rtimes_{\sigma,r}K}{Y}{C^\sigma}$ are continuous.  In
  particular, $C^{\tau}\rtimes_{\bar\sigma,r}K$ and $C^{\sigma}\rtimes
  _{\bar\tau,r}H$ are Morita equivalent.
\item[(2)] (Properness) The action
  $(\bar\sigma,\sigma,\sigma\rtimes\id)$ of $K$ on
  $\metrip{C^\tau}{X}{C\rtimes_{\tau,r}H}$ is proper with respect to
  $\metrip{(D_0,\bar\sigma)}{(X_0,\sigma)}{(E_0,\sigma\rtimes\id)}$ if
  for all $u,v,w,x\in X_0$
\begin{enumerate}
\item the function $t\mapsto \lip {C^\tau}<x,\sigma_t(v)>$ and its
  product with $\Delta_K(t)^{-1/2}$ are in $L^1(K, C^\tau)$;
\item the function $s\mapsto\Delta_H(s)^{-1/2}\rip
  C^\sigma<x,\tau_s(v)> $ is in $L^1(H,C^\sigma)$;
\item the integral $\int_H\int_K\tau_r(w)\sigma_t(\tau_r(x^*)v)\, dt\,
  dr$ is in $X_0$; and
\item the function $(r,s,t)\mapsto u\tau_r(v^*)\sigma_t\big(
  \tau_r(w^*)\tau_s(x) \big)\Delta_H(s)^{-1/2}$ is integrable.
\end{enumerate}
\item[(3)] (Tensor decomposition isomorphism) If in addition to (1)
  and (2) above
\begin{enumerate}
\item the action $\sigma\rtimes\id$ of $K$ on $C\rtimes_{\tau,r}H$ is
  saturated with respect to $E_0$; and
\item $E_0\subset L^1(H, Y)$
 \end{enumerate}
 then $X_0$ completes to a $C^\tau\rtimes_{\bar\sigma,r}K\sme
 C^\sigma\rtimes_{\bar\tau,r}H$-imprimitivity bimodule and
\[
(X\rtimes_{\sigma,r}K)\otimes (Y\rtimes_{\tau,r}H)\cong \overline{X_0}
\]
as $C^\tau\rtimes_{\bar\sigma,r}K\sme
C^\sigma\rtimes_{\bar\tau,r}H$-imprimitivity bimodules.
\end{enumerate}
\end{thm}
\begin{proof}
  Item~(1) follows from Proposition~\ref{prop-tensor-ib}, item~(2)
  from Proposition~\ref{prop-symmetric-proper}, and item~(3) from
  Proposition~\ref{prop-saturation} and \cite[Theorem~4.1]{hrw:xx}.
\end{proof}

\begin{cor}\label{cor-ibm}
  Suppose that $\tau:H\to\Aut C$ and $\sigma: K\to\Aut C$ are
  commuting actions which are proper and saturated with respect to the
  same dense $*$-subalgebra $C_0$, and suppose that the hypotheses of
  Theorem~\ref{thm-symmetric}(1)--(3) are satisfied.  Then $C_0$
  completes to a $C^\tau\rtimes_{\bar\sigma,r}K\sme
  C^\sigma\rtimes_{\bar\tau,r}H$-imprimitivity bimodule where the
  actions and inner products are given on dense objects by
  Equations~\ref{eq1-cor}--\ref{eq4-cor} below.
\end{cor}

\begin{proof}
  The action $(\bar\sigma,\sigma,\sigma\rtimes\id)$ of $K$ on
  $\metrip{C^\tau}{X}{C\rtimes_{\tau,r}H}$ is continuous by Item (1)
  and is proper with respect to
  $\metrip{(D_0,\bar\sigma)}{(C_0,\sigma)}{(E_0,\sigma\rtimes\id)}$ by
  Item (2).  The action $\sigma\rtimes\id$ of $K$ on
  $C\rtimes_{\tau,r}H$ is saturated with respect to $E_0$ by Item
  (3a), and since $E_0=\langle C_0\,,\,
  C_0\rangle_{C\rtimes_{\tau,r}H}$ it follows from
  \cite[Theorem~4.1]{hrw:xx} that the action $\sigma$ on $X$ is
  saturated with respect to $C_0$.  By \cite[Theorem~3.16]{hrw:xx} the
  completion of $C_0$ is a $C^\tau\rtimes_{\bar\sigma,r}K\sme
  (C\rtimes_{\bar\tau,r}H)^{\sigma\rtimes\id}$-imprimitivity bimodule.
  Finally, Items (2d), (3b), and Proposition~\ref{prop-saturation}
  allows us to identify $(C\rtimes_{\tau,r}H)^{\sigma\rtimes\id}$ and
  $C^\sigma\rtimes_{\bar\tau,r}H$.
  
  Chasing through the construction and identification above we can
  write down the actions and inner products for the
  $C^\tau\rtimes_{\bar\sigma,r}K\sme
  C^\sigma\rtimes_{\bar\tau,r}H$-imprimitivity bimodule obtained by
  competing $C_0$.  Let $x,y,c\in C_0$.  By \cite[Lemma~3.17]{hrw:xx},
\begin{align}
  {}_{C^\tau\rtimes_{\bar\sigma,r}K}\langle x\,,\,y\rangle\cdot c
  &=\int_K{}_{C^\tau\rtimes_{\bar\sigma,r}K}\langle
  x\,,\,y\rangle(t)\sigma_t(c)\Delta_K(t)^{-1/2}\, dt\\\notag &=\int_K
  {}_{C^\tau}\langle x\,,\,\sigma_t(y)\rangle\sigma_t(c)\, dt\\\notag
  &=\int_K\int_H \tau_s(x\sigma_t(y^*))\sigma_t(c)\, ds\,
  dt,\label{eq1-cor}
\end{align}
and now
\begin{align}
  c\cdot\langle x\,,\,y\rangle_{C^\sigma\rtimes_{\bar\tau,r}H}
  &={}_{C^\tau\rtimes_{\bar\sigma,r}K}\langle c\,,\,x\rangle\cdot
  y\\\notag &=\int_K\int_H \tau_s(c\sigma_t(x^*))\sigma_t(y)\, ds\,
  dt\label{eq2-cor}.
\end{align}
By \cite[Theorem~3.16]{hrw:xx} the left inner product is
\begin{equation}\label{eq3-cor}
{}_{C^\tau\rtimes_{\bar\sigma,r}K}\langle x\,,\,y\rangle(t)=\Delta_K(t)^{-1/2}{}_{C^\tau}\langle
x\,,\, \sigma_t(y)\rangle,
\end{equation}
and the right inner product is defined using the isomorphism of
$(C\rtimes_{\tau,r}H)^{\sigma\rtimes\id}$ onto
$C^\sigma\rtimes_{\bar\tau,r}H$ and Equation~\ref{eq-ripmult}:
\begin{equation}\label{eq4-cor}
\langle x\,,\, y\rangle_{C^\sigma\rtimes_{\bar\tau,r}H}(s)=\Delta_H(s)^{-1/2}\langle x\,,\,
\tau_s(y)\rangle_{C^\sigma}.
\end{equation}

\end{proof}

%%%%%%%%%%%%%%%%%%%%%%%%%%%%%%%%%%%%%%%%%%%%%%%%%%%%%%%%%%%%%%%%%%%%!!!

\section{Examples}\label{examples}

Our examples are based on the proper actions constructed in
\cite[Section~5]{rie:99}. There Rieffel starts with a proper action of
$G$ on the left of a locally compact Hausdorff space $P$, a nondegenerate
homomorphism $\theta:C_0(P)\to M(A)$, and an action $\alpha$ of $G$ on $A$ such that 
\[\alpha_s(\theta(f)a)=\theta(\lt_s(f))\alpha_s(a).
\]  
Rieffel proves in
\cite[Theorem~5.7]{rie:99} that $\alpha$ is proper in the sense of
\cite{rie:pm88} with respect to the subalgebra
\[
A_0:=\theta(C_c(P))A\theta(C_c(P))=\sp\set{\theta(f)a\theta(g):a\in
  A,\;f,g\in C_c(P)}.
\] 
The
homomorphism $\theta$ does \emph{not} necessarily have range in
$ZM(A)$, and consequently this setup includes striking examples (see Remark~\ref{rem-coactions}).  If $\theta$ is central then $A$ is a $C_0(P)$-algebra, and Proposition~\ref{prop-ex-tensor} and Theorem~\ref{thm-ex-symmetric} below reduce to results in \cite{hrw:tams00}.

Rieffel also says in \cite{rie:99} that the action is saturated if the action of $G$ on
$P$ is free; this is the content of the following lemma which we
prove in Appendix~\ref{app3}.

\begin{lemma}[Rieffel]\label{lem-dana-sat}
  Suppose $G$ acts freely and properly on a locally compact Hausdorff
  space $P$, there is a nondegenerate homomorphism $\theta:C_0(P)\to
  M(A)$, and an action $\alpha$ of $G$ on $A$ such that
  $\alpha_s(\theta(f)a)=\theta(\lt_s(f))\alpha_s(a)$.  Then the proper
  action $\alpha$ of $G$ on $A$ is saturated with respect to
%  $A_1:=\theta(C_c(P))A_0\theta(C_c(P))$, where
  $A_0:=\theta(C_c(P))A\theta(C_c(P))$.
\end{lemma}

%%%%%%%%%%%%

For our example, consider commuting proper actions of $K$ and $H$ on
the left and right of $P$. Suppose that we have a nondegenerate
homomorphism $\theta:C_0(P)\to M(C)$, and that we have commuting
actions $\sigma:K\to \Aut C$ and $\tau:H\to \Aut C$ satisfying
\begin{equation}\label{compatibility}
\bar\sigma_t(\theta(f))=\theta(\lt_t(f))
\quad\text{for $t\in K,$ and}\quad
 \bar\tau_s(\theta(f))=\theta(\rt_s(f))\mbox{ for $s\in H$}.
\end{equation} 
\cite[Theorem~5.7]{rie:99} implies that both $\sigma$ and $\tau$ are
proper with respect to the same subalgebra
$C_0:=\theta(C_c(P))C\theta(C_c(P))$.

\begin{prop}\label{prop-ex-tensor}
Suppose that $\sigma:K\to \Aut C$ and $\tau:H\to \Aut C$ are commuting actions and that the actions of $K$ and $H$ on $P$ are free and proper.  If $C_0(P)$ maps bi-equivariantly into $M(C)$ then $C^\tau\rtimes_{\bar\sigma,r} K$
  and $C^\sigma\rtimes_{\bar\tau,r} H$ are Morita equivalent.
\end{prop}

By Theorem~\ref{thm-symmetric}(1), it suffices to verify that the
actions $(\bar\sigma,\sigma,\sigma\times\id)$ and
$(\tau\rtimes\id,\tau,\bar\tau)$ of $K$ and $H$ on
$\metrip{C^\tau}X{C\rtimes_{\tau,r}H}$ and
$\metrip{C\rtimes_{\sigma,r}K}Y{C^\sigma}$ are continuous.
Furthermore, if $B:= C\rtimes_{\sigma\times\tau, r}(H\times K)$ then
$\big(X\rtimes_{\sigma,r}K\big)\otimes_B\big( Y\rtimes_{\tau,r}H\big)$
is a $C^\tau\rtimes_{\bar\sigma,r}K\sme
C^\sigma\rtimes_{\bar\tau,r}H$-imprimitivity bimodule by
Proposition~\ref{prop-tensor-ib}.  We retain the notation from
Sections~\ref{tensor-result} and~\ref{concrete-result}, and we will drop all mention of
$\theta$ --- we must remember that $fc\neq cf$.  The key to our
calculations is the following lemma.

\begin{lemma}\label{cpctsupp}
  For $x=fbg$ and $y=hck$ in $C_0$, the functions $s\mapsto
  x\tau_s(y)$ and $t\mapsto x\sigma_t(y)$ have compact support
  depending only on $\supp g$ and $\supp h$.
\end{lemma}

\begin{proof}
  The consistency conditions (\ref{compatibility}) imply that
  $x\tau_s(y)=fbg\bigl(\rt_s(h)\bigr)\tau_s(ck)$.  But $g\rt_s(h)$ is
  nonzero only if
\[
\supp g\cap \supp(\rt_s(h))=\supp g\cap (\supp h)s^{-1}
\]
is nonempty.  Therefore the support of $s\mapsto x\tau_{s}(y)$ is
contained in
\[
\{s\in H: \supp g\cap (\supp h)s^{-1}\not=\emptyset\}
\]
which is compact because $H$ acts properly on $P$.  The other part is
similar.
\end{proof}

\begin{proof}
  [Proof of Proposition~\ref{prop-ex-tensor}] We want to show that
  $(\bar\sigma,\sigma,\sigma\rtimes\id)$ is a continuous action on
  $\metrip{C^\tau}X{C\rtimes_{\tau,r}H}$, and
  Theorem~\ref{thm-symmetric}(1) implies it suffices to see that for
  each fixed $x\in X_0=C_0$ the map $t\mapsto \brip C\rtimes_{\tau,r}H
  <\sigma_t(x),x>$ is continuous.  Note that
\begin{align}
  \|\brip \cctrh<\sigma_t(x),x>- \rip \cctrh <x, x>\|
  &\leq\int_H \|\brip \cctrh<\sigma_t(x)-x, x>(s)\|\, ds\notag\\
  &=\int_H \|(\sigma_t(x)-x)^*\tau_s(x)\|\Delta_H(s)^{-1/2}\,
  ds\label{eq-cty}.
\end{align}
We claim that the integrand in \eqref{eq-cty} has support in a compact
set $L$ whenever $t$ is in a sufficiently small neighborhood $M$ of
$e$, so that
\[
\eqref{eq-cty}\leq\|\sigma_t(x)-x \|
\|\tau_s(x)\Delta_H(s)^{-1/2}\|_\infty\mu_H(L)\to 0\text{\ as\ }t\to 0
\]
because $\sigma:K\to\Aut C$ and $s\mapsto \tau_s(x)\Delta_H(s)^{-1/2}$
are continuous.

To prove our claim, set $x=fcg$, and choose $h,k\in C_c(P)$ such that
$h=1$ on a neighborhood of $\supp f$ and $k=1$ on a neighborhood of
$\supp g$.  Then there exists a neighborhood $M$ of $e$ in $K$ such
that it $t\in M$ then $h=1$ on $\supp\lt_t(f)$ and $k=1$ on
$\supp\lt_t(g)$, so that
\[
\sigma_t(x)-x=h(\sigma_t(x)-x)k
\]
Now, by Lemma~\ref{cpctsupp}, when $t\in M$ the support of
$(\sigma_t(x)-x)^*\tau_s(x)$ is compact and depends only on $k$ and
$f$. This proves the claim.  Note that by the symmetry of our
situation the action $(\tau\rtimes\id,\tau,\bar\tau)$ is continuous on
$\metrip{C\rtimes_{\sigma,r}K}Y{C^\sigma}$.
\end{proof}

In fact, our Theorem~\ref{thm-symmetric} gives us two
$C^\tau\rtimes_{\bar\sigma,r} K$--$C^\sigma\rtimes_{\bar\tau,r}
H$-imprimitivity bimodules and an isomorphism between the two:

\begin{thm}\label{thm-ex-symmetric} 
Suppose that $\sigma:K\to \Aut C$ and $\tau:H\to \Aut C$ are commuting actions and that the actions of $K$ and $H$ on $P$ are free and proper.  If $C_0(P)$ maps bi-equivariantly into $M(C)$
then all the hypotheses of
  Theorem~\ref{thm-symmetric} are satisfied.  That is, the two
  $C^\tau\rtimes_{\bar\sigma,r} K \sme
  C^\sigma\rtimes_{\bar\tau,r}H$-imprimitivity bimodules
  $\big(X\rtimes_{\sigma,r}K\big)\otimes_B\big(
  Y\rtimes_{\tau,r}H\big)$ and $\overline{ X_0}$ are isomorphic.
\end{thm}

\begin{proof} In view of the proof of Proposition~\ref{prop-ex-tensor} 
  we only need to verify Items 2~and 3 of Theorem~\ref{thm-symmetric}.
 
  2. (Properness) To see that $(\bar\sigma,\sigma,\sigma\rtimes\id)$
  is a proper action on $\metrip {C^\tau}X{C\rtimes_{\tau,r}H}$ with
  respect to $\metrip{D_0}{(X_0)}{E_0}$ we need to check items
  (a)--(d) of Theorem~\ref{thm-symmetric}(2).  For parts (a)~and (b),
  note that if $x,v\in X_0$, then $t\mapsto x\sigma_t(v)$ has compact
  support by Lemma~\ref{cpctsupp}. Since
\[
\blip {C^\tau}<x,\sigma_t(v)>c=\int_H\tau_s(x\sigma_t(v^{*}))c\, ds
\] 
for every $c\in C_0$, it follows that the continuous function
$t\mapsto \blip {C^\tau}<x,\sigma_t(v)>$ has compact support with norm
at most $\|x\|\|v\|$.  Thus it and its product with
$\Delta_K(t)^{-1/2}$ are integrable. By the symmetry of our situation
the function $s\mapsto\Delta_H(s)^{-1/2}\brip C^\sigma<x, \tau_s(v)>$
is in $L^1(H, C^\sigma)$.

For (c), we need
%$\rip (C\rtimes_{\tau,r}H)^{\sigma\rtimes\id}<x,
%v>$ to multiply $E_0$
%and $x\cdot\rip (C\rtimes_{\tau,r}H)^{\sigma\rtimes\id}<v,
%w>$ to be back in
%$X_0$, and it suffices to check that
\begin{equation}\label{eq-multiplier}
\int_H\int_K\tau_r(w)\sigma_t(\tau_r(x^*)v)\, dt\, dr\in X_0  \quad 
\text{whenever $v,w,x\in X_0$.}
\end{equation}
Note that $r\mapsto \tau_r(x^*)v$ has compact support $L$, say, by
Lemma~\ref{cpctsupp}. Let $w=fcg$ and choose $h\in C_c(P)$ such that
$h=1$ on $(\supp f)\cdot L^{-1}$. Then $\tau_r(w)=h\tau_r(w)$ for
$r\in L$ so that the left hand side of \eqref{eq-multiplier} equals
\[
h\int_H\int_K\tau_r(w)\sigma_t(\tau_r(x^*)v)\, dt\, dr,
\]
and is therefore an element of $C_c(P) C$.  Since $t\mapsto
w\sigma_t(x^*)$ is also compactly supported we just repeat the
argument for the right side of the integral to get that
$\int_H\int_K\tau_r(w)\sigma_t(\tau_r(x^*)v)\, dt\, dr\in X_0=
C_{c}(P) C C_{c}(P)$.

For part~(d), we check that
\begin{equation*}
(r,s,t)\mapsto u\tau_r(v^*)\sigma_t\big( \tau_r(w^*)\tau_s(x)
\big)\Delta_H(s)^{-1/2}
\end{equation*}
is integrable whenever $u,v,w,x\in X_0$.  Again, by
Lemma~\ref{cpctsupp}, the maps $r\mapsto u\tau_r(v^*)$ and $t\mapsto
v^*\sigma_t(w^*)$ have compact supports.  Thus it suffices to see that
for $r$ in a compact set $L$ the function $s\mapsto
\tau_r(w^*)\tau_s(x)$ has compact support. If $w^*=fcg$ and $x=kdh$
where $f,g,k,h\in C_c(P)$, then
$\tau_r(w^*)\tau_s(x)=\tau_r(fc)\tau_r(g)\tau_s(k)\tau_s(dh)$ and
$\tau_r(g)\tau_s(k)$ is nonzero if and only if $(\supp g)\cdot
L^{-1}\cap (\supp k)\cdot s^{-1}\neq\emptyset$. Since the action of
$H$ on $P$ is proper the set $\{s\in H:(\supp g)\cdot L^{-1}\cap
(\supp k)\cdot s^{-1}\neq\emptyset\}$ is compact.

3. (Tensor decomposition isomorphism) That the action of
$\sigma\rtimes\id$ of $K$ on $C\rtimes_{\tau,r}H$ is saturated with
respect to $E_0$ follows by noting that $E_0=C_c(P)E_0C_c(P)$ and
applying Lemma~\ref{lem-dana-sat} with $A_0=E_0$. To see that
$E_0=C_c(P)E_0C_c(P)$, let $x=fcg$ and $y=kdh\in C_0$ and let $L$ be
the compact support of $x\tau_s(y)$. Choose $l\in C_c(P)$ such that
$l$ is identically one on $(\supp h)\cdot L^{-1}$ and on $\supp f$.
Then $x\tau_s(y)=fcg\tau_s(kdh)=l\big(fcg\tau_s(kdh)\big)l$.

Finally, if $e\in E_0$ is given by
$e(s)=\Delta_H(s)^{-1/2}x\tau_s(y^*)$, then
\begin{align}\label{eq-in-l-one}
  \|e\|_{L^1(H,Y)} &=\int_H\|e(s)\|_Y\, ds =\int_H\| \blip
  {C\rtimes_{\sigma,r}K}
  <e(s), e(s)>\|^{1/2}\, ds\notag\\
  &\leq\int_H\int_K\| \blip {C\rtimes_{\sigma,r}K}<e(s), e(s)>
  (t)\|_C^{1/2}\, dt\,
  ds\notag\\
  &=\int_H\int_K\| e(s)\sigma_t(e(s)^*)\|_C^{1/2} \Delta_K(t)^{-1/4}\,
  dt\, ds\notag\\
  &=\int_H\int_K\|
  x\tau_s(y^*)\sigma_t(\tau_s(y)x^*)\|_C^{1/2}\Delta_H(s)^{-1/2}
  \Delta_K(t)^{-1/4}\, dtds <\infty
\end{align}
because the integrand is continuous with compact support (because
$s\mapsto x\tau_s(y^*)$ and $t\mapsto y^*\sigma_t(y)$ have compact
supports).  Hence $E_{0}\subset L^{1}(H,Y)$ as required.
\end{proof}

\begin{remark} \label{rem-coactions}
  It has apparently not been noticed that Rieffel's construction in
  \cite[Theorem~5.7]{rie:99} implies that the dual action on any
  crossed product by a coaction is proper.
To see this, suppose $\delta:A\to M(A\otimes C^*(G))$ is a coaction of a locally compact group $G$ on a $C^*$-algebra $A$.  Then the crossed product $A\rtimes_\delta G$ is generated by a universal covariant representation $(j_A, j_{C(G)})$ of $(A, C_0(G))$ in $M(A\rtimes_\delta G)$.  Since the dual action $\hat\delta:G\to \Aut(A\rtimes_\delta G)$ is characterized by
\[
\hat\delta_s\bigl(j_A(a)j_{C(G)}(f)\bigr)=j_A(a)j_{C(G)}(\rt_s(f)),
\]
the homomorphism $j_{C(G)}$ is equivariant for the actions $\hat\delta$ and $\rt:G\to\Aut(C_0(G))$.
Applying \cite[Theorem~5.7]{rie:99} to $j_{C(G)}$ shows that $\hat\delta$ is a proper action.  More generally, it shows that $\hat\delta|_H$ is proper for any closed subgroup 
$H$ of $G$; this improves a result of Mansfield \cite[Theorem~30]{man:jfa91} for normal amenable $H$.

We can therefore apply Rieffel's original theorem from \cite{rie:pm88} to obtain a Morita equivalence between $(A\rtimes_\delta G)\rtimes_{\hat\delta, r}H$ and a generalized fixed-point algebra $(A\rtimes_\delta G)^H$. 
Since $A\rtimes_\delta G$ is generated by the
universal covariant representation $(j_A,j_{C(G)})$, and is even spanned by elements
of the form $j_A(a)j_{C(G)}(f)$ (see \cite[\S2]{rae:plms92}), it is
tempting to guess that $(A\rtimes_\delta G)^H$ is at least generated
by elements of the form $j_A(a)j_{C(G)}(f)$ for $f\in C_0(G/H)$, and
hence coincides with the candidate for the crossed product
$A\rtimes_\delta G/H$ by the homogeneous space discussed in
\cite{ekr:jot98}.  This is indeed the case if $G$ is discrete
\cite{dpr}.  Thus, Rieffel's theorem could
give an extension of Mansfield's imprimitivity theorem to coactions of
arbitrary homogeneous spaces (as opposed to quotients by normal
amenable subgroups as in \cite{man:jfa91}; see \cite{kalqui:mpcps98}
for a discussion of this problem for non-amenable normal subgroups).
Unfortunately it does not seem to be easy to write a typical element
of $C_c(G)(A\rtimes_\delta G)C_c(G)$ in the form $j_A(a)j_{C(G)}(f)$,
or to do so approximately in such a way that one can verify the
existence of the multiplier $\rip D<x, y>$. Indeed, it is the content
of one of Mansfield's main theorems \cite[Theorem 19]{man:jfa91}, that
there are such multipliers when $H$ is normal and amenable and $x,y$
lie in a dense subalgebra ${\mathcal{D}}$ of $A\rtimes_\delta G$, and
the proof of this theorem relies on some very subtle estimates.
This analysis is a crucial ingredient in the proof of \cite[Theorem 30]{man:jfa91}.

Pask and Raeburn showed in \cite{PR} that a free action on a directed
graph $E$ induces a proper action on the associated graph algebra
$C^*(E)$. In their result too there is
an underlying proper $G$-space $P$ together with a non-central equivariant
map of $C_0(P)$ into $M(C^*(E))$: just take
$P$ to be the set of vertices of the graph with the discrete topology.  

Thus all the main examples of proper actions come with the existence of an underlying proper action on a space.
\end{remark}

\appendix

\section{Proof of Proposition~\ref{prop-mult}}\label{appendix1}

The object of these appendices is to make sense of certain
manipulations with vector-valued integrals needed to give careful
proofs of Proposition~\ref{prop-mult} and Lemma~\ref{lem-variable}.
If $A$ is a \cs-algebra, then the collection of Bochner-integrable
functions from $G$ to $A$ will be denoted by $\L^{1}(G,A)$, and the
Banach space of equivalence classes of integrable functions agreeing
almost everywhere will be denoted by $L^{1}(G,A)$.

For motivation, recall that we can realize $A\rtimes_{\alpha}G$ as the
enveloping $C^*$-algebra of the Banach $*$-algebra $L^{1}(G,A)$.  The
product and involution are given by the usual formulas:
\begin{gather}\label{eq-a1}
  f*g(s):=\int_{G}f(r)\alpha_{r}\bigl(g(r^{-1}s)\bigr)\, ds, \text{
    and}\\ \notag f^{*}(s):=
  \Delta(s^{-1})\alpha_{s}\bigl(f(s^{-1})^{*}\bigr)\quad\text{for
    $f,g\in \mathcal{L}^{1}(G,A)$.}
\end{gather}
It takes some work, though, to see that $f^{*}$ and $f*g$ are
well-defined elements of $L^{1}(G,A)$.  The first step is to
see that $(r,s)\mapsto f(r)\alpha_{r}\bigl(g(r^{-1}s)\bigr)$ is a
measurable function from $G\times G$ to $A$.  This is a bit thorny as
there is no \emph{a priori} reason to suspect that $(r,s)\mapsto
g(r^{-1}s)$ is measurable if $g$ is merely measurable rather than
continuous or Borel.  There are a number of finesses for this.  Here,
we use the following lemma; we assume $g$ is integrable to
ensure that we can approximate it globally with functions in
$C_{c}(G,A)$.

\begin{lemma}
  \label{lem-single}
  Let $(A,G,\alpha)$ be a dynamical system.  Suppose that
  $g\in\L^{1}(G,A)$ and $h(r,s):=\alpha_{r}\bigl(g(r^{-1}s)\bigr)$.
  Then $h:G\times G\to A$ is measurable.
\end{lemma}
\begin{proof}
  Since Haar measure is a Radon measure, a subset $S$ of $G$ is
  measurable if and only if $S\cap C$ is measurable for all compacts
  sets $C\subset G$ \cite[Theorem~III.11.31]{hr:abstract}. It follows
  that $h$ is measurable if and only if $h|_{L}$ is measurable for
  each compact set $L\subset G\times G$.  Therefore it will suffice to
  show that $h|_{K\times K}$ is measurable for each compact set
  $K\subset G$.  To do this, we'll produce measurable functions
  $h_{n}$ such that $h_{n}\to h$ almost everywhere on $K\times K$.
  
  Since $g\in\L^{1}(G,A)$, there are $g_{n}\in C_{c}(G,A)$ such that
  $g_{n}\to g$ in $L^{1}(G,A)$.  Passing to a subsequence and
  relabeling, we can assume that there is a \emph{Borel} null set $N$
  such that $g_{n}(s)\to g(s)$ for all $s\notin N$. Since $g_{n}$ is
  continuous, $h_{n}(r,s):=\alpha_{r}\bigl(g_{n}(r^{-1}s)\bigr)$
  defines a measurable function (continuous in fact), and
  \begin{equation*}
    h_{n}(r,s)\to h(r,s)
  \end{equation*}
  for all $(r,s)\in K\times K\setminus D$, where $D=\{(r,s)\in K\times
  K:r^{-1}s\in N\}$.  Since $N$ is Borel, $D$ is a measurable subset
  of $K\times K$.  Since $K\times K$ has finite product-measure,
  Tonelli's Theorem, as proved in
  \cite[Theorem~III.13.9]{hr:abstract}, implies that
\begin{equation*}
  \mu\times \mu(D) =\int_{G} \mu(D_{r})\, dr,
\end{equation*}
where $D_{r}:=\{s:(r,s)\in D\}$.  Since $D_{r}\subset rN$ and
$\mu(rN)=0$ for all $r$, it follows that $D$ is a null set.  This
completes the proof.
\end{proof}

With Lemma~\ref{lem-single} in hand, the measurability of the function
$m$ given by $(r,s)\mapsto f(r)\alpha_{r}\bigl(g(r^{-1}s)\bigr)$
follows because the product of vector-valued measurable functions is
measurable.\footnote{It suffices, for example, to see that product of
  measurable simple functions is again a measurable simple function.}
Since $f$ and $g$ are integrable, $m$ must be supported on a
$\sigma$-finite set, and Tonelli's Theorem shows that $\|m\|$ is in
$\L^{1}(G\times G)$.  Consequently, $m\in \L^{1}(G\times G,A)$.  Now a
vector-valued Fubini's Theorem, such as
\cite[Theorem~II.16.3]{fd:representations}, implies that the
right-hand side of \eqref{eq-a1} is defined for almost all $s$ and
that $f*g$ is a well-defined element of $L^{1}(G,A)$.  (Since
$\|f*g\|_{1}\le \|f\|_{1}\| g\|_{1}$, it is not hard to see that the
class of $f*g$ depends only on the classes of $f$ and $g$.)

\begin{remark}
  \label{rem-three-fcns}
  Similar considerations are often glossed over when it is observed
  that convolution is associative.  For example, if $f$, $g$ and $h$
  are in $\L^{1}(G,A)$, then Lemma~\ref{lem-single} implies that
  \begin{equation}\label{eq-a5}
    (r,t,s)\mapsto
    f(r)\alpha_{r}\bigl(g(r^{-1}t)\bigr)\alpha_{t}\bigl(
    h(t^{-1}s)\bigr) 
  \end{equation}
  is measurable, and Tonelli's Theorem implies that \eqref{eq-a5} is
  integrable on $G\times G\times G$.  Then Fubini's Theorem implies
  that for almost all $s\in G$,
\begin{equation}
  \label{eq-a6}
  (r,t)\mapsto f(r)\alpha_{r}\bigl(g(r^{-1}t)\bigr)\alpha_{t}\bigl(
    h(t^{-1}s)\bigr) 
\end{equation}
is in $\L^{1}(G\times G,A)$.  This will allow us to apply Fubini's
Theorem to double integrals with integrands such as \eqref{eq-a6} in
the sequel.\footnote{It is possible (by approximating by simple
  functions) to see that \eqref{eq-a6} is measurable without resorting
  to functions on $G\times G\times G$.  However, it is interesting to
  note that it is not obvious that \eqref{eq-a6} is integrable without
  appealing to the integrability of \eqref{eq-a5}.}
\end{remark}

We can view the multiplier algebra $M(A)$ as the $C^*$-algebra
${\mathcal L}(A_{A})$ of adjointable operators on the right Hilbert
$A$-module $A_A$.

\begin{lemma}
  \label{lem-l-one-mult}
  Suppose that $(A,G,\alpha)$ is a dynamical system and that $T$ and
  $S$ are bounded linear operators on $L^{1}(G,A)$ such that for all
  $f$ and $h$ in $L^{1}(G,A)$ we have
  \begin{equation*}
    T(f*h)=Tf*h,\quad S(f*h)=Sf*h,\quad\text{and}\quad (Tf)^{*}*h
    = f^{*}*Sh.
  \end{equation*}
  Then $T$ and $S$ extend to elements of
  $\L(A\rtimes_{\alpha}G)=M(A\rtimes_{\alpha}G)$ satisfying $T^{*}=S$.
\end{lemma}
\begin{proof}
  Let $\{e_{i}\}$ be a bounded approximate identity for $L^{1}(G,A)$.
  Then if $\pi$ is a representation of $A\rtimes_{\alpha}G$,
  \begin{align*}
    \bigl\|\pi(Tf)\bigr\|&=\lim_{i}\bigl \|\pi\bigl(T(e_{i}*f)\bigr) \bigr\| \\
    &\le \limsup_{i} \bigl\|\pi\bigl(Te_{i})\bigr)\bigr\|\bigl\|\pi(f)\bigr\| \\
    &\le M\|T\| \bigl\|\pi(f)\bigr\|.
  \end{align*}
  It follows that $T$ is bounded with respect to the universal norm on
  $L^{1}(G,A)\subset A\rtimes_{\alpha}G$.  Thus $T$ and $S$ extend to
  operators on $A\rtimes_{\alpha}G$.  Since
\begin{equation*}
 \rip A\rtimes_\alpha G<Tf, h> = (Tf)^*h=f^*Sh=
\rip A\rtimes_\alpha G<f, Sh>
\end{equation*}
it follows that $T$ is adjointable with $T^{*}=S$.
\end{proof}

\begin{remark}
  \label{rem-m-convolutions} 
  If $B$ is any $C^*$-subalgebra of $M(A)$ then we can identify
  $\L^{1}(G,B)$ with a subalgebra of $L^{1}(G,M(A))$.  In particular,
  if $f\in \L^{1}(G,A)$ and $g\in \L^{1}(G,B)$, then
  Lemma~\ref{lem-single} implies that
  \begin{equation*}
    (r,s)\mapsto f(r)\bar\alpha_{r}\bigl(g(r^{-1}s)\bigr)
  \end{equation*}
  is a measurable function of $G\times G$ into $M(A)$ taking values in
  $A$.  Thus it is a measurable function of $G\times G$ into $A$.  A
  similar statement can be made if $f\in \L^{1}(G,B)$ and $g\in
  \L^{1}(G,A)$.
\end{remark}

\begin{proof}
  [Proof of Proposition~\ref{prop-mult}] The integrand on the
  right-hand side of \eqref{eq:5} is measurable in view of
  Remark~\ref{rem-m-convolutions} and Lemma~\ref{lem-single}.  Thus
  applications of the Tonelli and Fubini Theorems imply that the
  right-hand side of \eqref{eq:5} defines an element $T_{g}f$ in
  $L^{1}(G,A)$.  Furthermore, $\|T_{g}\|\le\|g\|_{1}$.
  
  If $h\in\mathcal{L}^{1}(G,A)$, then $f*h$ is too, and by definition,
  for almost all $s$,
\begin{align*}
  T_{g}(f*h)(s) &= \int_{G}g(r)\alpha_{r}\bigl(f*h(r^{-1}s)\bigr)\, dr
  \\
  &= \int_{G}\int_{G} g(r)\alpha_{r}\bigl(f(t)\alpha_{t}\bigl(
  h(t^{-1}r^{-1 }s)\bigr)\bigr)\, dt\, dr \\
  &= \int_{G}\int_{G} g(r)\alpha_{r} \bigl(f(r^{-1}t)\bigr)
  \alpha_{t}\bigl( h(t^{-1}s)\bigr) \, dt\, dr.
\end{align*}
The integrand $(r,t)\mapsto g(r)\alpha_{r}\bigl( f(r^{-1}t)\bigr)
\alpha_{t}\bigl(h(t^{-1}s)\bigr) $ is in in $\L^{1}(G\times G,A)$ by
Remark~\ref{rem-three-fcns}, so Fubini's Theorem implies that for
almost all $s$,
\begin{align*}
  T_{g}(f*h)(s)&= \int_{G} T_{g}f(t)\alpha_{t}\bigl(h(t^{-1}s)\bigr) \, dt \\
  &= T_{g}f*h(s).
\end{align*}
Thus $T_{g}(f*h) = T_{g}f*h$ in $L^{1}(G,A)$.

Next we want to show that $(T_{g}f)^{*}*h = f^{*}*(T_{g^{*}}h)$ in
$L^{1}(G,A)$.  But for almost all $s$,
\begin{align}
  (T_{g}f)^{*}*h(s) &= \int_{G}(T_{g}f)^{*}(r)
  \alpha_{r}\bigl(h(r^{-1}s)\bigr) \, dr\notag \\
  &= \int_{G} \alpha_{r}\bigl(T_{g}f(r^{-1})\bigr)^{*}\Delta(r^{-1})
  \alpha_{r} \bigl(h(r^{-1}s)\bigr) \, dr\notag \\
  &= \int_{G}\int_{G} \alpha_{r}\bigl(g(t)
  \alpha_{t}\bigl(f(t^{-1}r^{-1})\bigr)\bigr)^{*}\Delta(r^{-1})
  \alpha_{r} \bigl(h(r^{-1}s)\bigr) \, dt\, dr \notag \\
  &= \int_{G} \int_{G} \alpha_{r^{-1}t}\bigl(f(t^{-1}r)^{*}\bigr)
  \alpha_{r^{-1}}\bigl(g(t)^{*}\bigr) \alpha_{r^{-1}}
  \bigl(h(rs) \bigr) \, dt\, dr. \notag \\
  \intertext{It follows from simple variations on
    Lemma~\ref{lem-single} that the integrand above is measurable, and
    the (scalar version of) Tonelli's Theorem implies it is
    integrable. Hence we can use the vector-valued version of Fubini's
    Theorem to conclude that, for almost all $s$,} (T_{g}f)^{*}*h(s)&=
  \int_{G}\int_{G} \alpha_{r^{-1}t} \bigl(f(t^{-1}r)^{*}\bigr)
  \alpha_{r^{-1}} \bigl(g(t)^{*}\bigr) \alpha_{r^{-1}}
  \bigl(h(rs)\bigr) \, dr\, dt \notag \\
  \intertext{and, since it now makes sense to send $r\mapsto tr$,
    this} &= \int_{G}\int_{G} \alpha_{r^{-1}}\bigl(f(r)^{*}\bigr)
  \alpha_{r^{-1}t^{-1}} \bigl(g(t)^{*}h(trs)\bigr) \, dr\, dt \notag
  \\
  \intertext{which, after sending $r\mapsto r^{-1}$ and $t\mapsto
    t^{-1}$, is} &= \int_{G}\int_{G} f^{*}(r)\alpha_{r} \bigl(
  g^{*}(t) \alpha_{t}
  \bigl(h(t^{-1}r^{-1}s)\bigr)\bigr) \, dr\, dt\notag. \\
  \intertext{To apply Fubini's Theorem we need to see that the above
    integrand is a measurable function of $r$, $t$ and $s$; as before,
    this follows from variations on Lemma~\ref{lem-single}.  Thus for
    almost all $s$,} (T_{g}f)^{*}*h(s) &= \int_{G}f^{*}(r) \alpha_{r}
  \Bigl( \int_{G} g^{*}(t) \alpha_{t}
  \bigl( h(t^{-1}r^{-1}s)\bigr) \, dt\Bigr) \, dr \label{eq:7} \\
  &= \int_{G} f^{*}(r) \alpha_{r}\bigl(T_{g^{*}}h(r^{-1}s)\bigr) \, dr
  \notag \\
  &= f^{*}*(T_{g^{*}}h)(s).\notag
\end{align}
It follows from Lemma~\ref{lem-l-one-mult} that $T_{g}$ defines a
multiplier with adjoint $T_{g}^{*}=T_{g^{*}}$.

Next we want to establish \eqref{eq-prop-mult}.  Fubini's Theorem
implies that the right-hand side of \eqref{eq-prop-mult} defines a
function $l$ in $\L^{1}(G,A)$.  Since $fT_{g}=(T_{g^{*}}f^{*})^{*}$ in
$A\rtimes_{\alpha}G$, it follows that $fT_{g}\in L^{1}(G,A)$, and we
have to show that $l=fT_{g}$ in $L^{1}(G,A)$.

Let $h\in\L^{1}(G,A)$.  We can repeat the computation of
$(T_{g}f)^{*}*h$ above with $f$ replaced by $f^{*}$ and $g$ replaced
by $g^{*}$, and use \eqref{eq:7} to conclude that
\begin{align*}
  (fT_{g})*h(s)=\bigl(T_{g^{*}}f^{*}\bigr)^{*}*h(s) &=
  \int_{G}\int_{G} f(r)\alpha_{r}\bigl(g(t)\alpha_{t}
  \bigl( h(t^{-1}r^{-1}s)\bigr)\bigr) \, dt\, dr \\
  &= \int_{G}\int_{G} f(r) \bar\alpha_{r}\bigl(g(r^{-1}t)\bigr)
  \alpha_{t}\bigl( h(t^{-1}s)\bigr) \, dt\, dr \\
  \intertext{which by Fubini's Theorem is} &= \int_{G}\Bigl( \int_{G}
  f(r) \bar\alpha_{r}\bigl(g(r^{-1}t)\bigr)
  \, dr\Bigr) \alpha_{t}\bigl(h(t^{-1}s)\bigr) \, dt \\
  &= \int_{G}l(t)\alpha_{t}\bigl(h(t^{-1}s)\bigr) \, dt \\
  &= l*h(s).
\end{align*}
It follows that $l=fT_{g}$ in $L^{1}(G,A)$, as claimed.
\end{proof}

\section{The Proof of 
  Lemma~\ref{lem-variable}}
\label{appendix2}

To begin with, let $C$ be an arbitrary $C^*$-algebra.  We let $H$ and
$K$ be arbitrary locally compact groups with Haar measures $\mu_{H}$
and $\mu_{K}$, respectively.  (For most of what follows, $C$ could be
any Banach space and $H$ and $K$ could be arbitrary locally compact
spaces equipped with Radon measures $\mu_{H}$ and $\mu_{K}$,
respectively.)

% The collection of integrable functions from $H$ to $C$ will be denoted
% by $\L^{1}(H,C)$; and the Banach space of equivalence classes of
% integrable functions agreeing almost everywhere will be denoted by
% $L^{1}(H,C)$. 
If $f\in L^{1}(H,C)$, then we will sometimes write
\begin{equation*}
  \intc_{H} f(s)\, ds
\end{equation*}
to emphasize where our integral takes its value.
% For example, if $h\in\slohc$ and $g\in
% \L^{1}\bigl(K,\lohc\bigr)$, then when we write
% \begin{equation*}
%   \intlohc_{K}g(t)\, dt(s)=h(s),
% \end{equation*}
% we are merely indicating that $h$ is a representative for the class of
% \begin{equation*}
%   \intlohc_{K}g(t)\, dt
% \end{equation*}
% in $\lohc$.

The key lemma is a very special case of
\cite[Lemma~III.11.17]{ds:linear}, and, modulo some facts about
vector-valued integrals, has a fairly straightforward proof.  We
include the proof here because the arguments from \cite{ds:linear} are
difficult and can be substantially simplified in our
situation.

\begin{lemma}
  \label{lem-dunford}
  Suppose that $g\in\L^{1}\bigl(K,\lohc\bigr)$ and that $f$ is a
  $\mu_{K}\times \mu_{H}$-integrable $C$-valued function on $K\times
  H$ such that for almost all $t$, the class of $f(t,\cdot)$ equals
  $g(t)$.  Then
\begin{enumerate}
\item for almost all $s\in H$, $f(\cdot,s)\in \L^{1}(K,C)$,
\item the function
\begin{equation*}
s\mapsto \intc_{K} f(t,s)\, dt
\end{equation*}
is in $\slohc$, and
\item as elements of $\lohc$,
  \begin{equation*}
    \intlohc_{K}g(t)\, dt=\Bigl(s\mapsto \intc_{K}f(t,s)\, dt\Bigr);
  \end{equation*}
  that is, for almost all $s$,
\begin{equation*}
  \intlohc_{K}g(t)\, dt(s)=\intc_{K}f(s,t)\, dt.
\end{equation*}
\end{enumerate}
\end{lemma}

The first two assertions follow immediately from any vector-valued
Fubini Theorem worthy of the name.  Our proof of the third assertion
is straightforward except for the following result which is
\cite[Lemma~III.6.8]{ds:linear}.

\begin{lemma}
  \label{lem-dsIII-6-8}
  If $f\in\slohc$ and if
  \begin{equation*}
    \int_{E}f(s)\, ds=0
  \end{equation*}
  for all measurable subsets $E\subset H$, then $f$ vanishes almost
  everywhere.
\end{lemma}

Although the proof of Lemma~\ref{lem-dsIII-6-8} is routine in the
scalar-valued case, we see no elementary proof in the vector-valued
case.  The proof in \cite[III \S2]{ds:linear} goes as follows.  Given
$h\in\slohc$, define a $C$-valued set function on measurable subsets
of $H$ by
\begin{equation*}
  \lambda(E):=\intc_{E}h(s)\, ds.
\end{equation*}
Of course, $\lambda$ is additive, and has a \emph{total variation}
defined by
\begin{equation*}
  \nu(E):=\sup\set{\sum_{i}\|\lambda(E_{i})\|:\text{$E_{1}$,\dots,$E_{n}$
  is a partition of $E$}}.
\end{equation*}
One sees easily that
\begin{equation*}
  \nu(E)\le \int_{E}\|h(s)\|\, ds.
\end{equation*}
But it can be shown \cite[Theorem~III.2.20]{ds:linear} that
\begin{equation*}%\label{eq-b3}
  \nu(E)=\int_{E}\|h(s)\|\, ds.
\end{equation*}
With these assertions in place, Lemma~\ref{lem-dsIII-6-8} is an easy
consequence; if $h=f$ as in Lemma~\ref{lem-dsIII-6-8}, then $\nu(E)=0$
for all $E$, and so $f$ is zero almost everywhere.

The only other tool we need for the proof of Lemma~\ref{lem-dunford}
is that bounded linear maps commute with vector-valued integrals.  For
each measurable subset $E\subset H$, we can define a bounded linear
map $\phi_{E}:\lohc\to C$ by
\begin{equation*}
\phi_{E}(h):=\intc_{E}h(s)\, ds
% :=\int_{H}^{C}\mathbb{I}_{E}(s)h(s)\, ds.  
\end{equation*}

\begin{proof}
  [Proof of Lemma~\ref{lem-dunford}] To prove the final assertion, it
  suffices, in view of Lemma~\ref{lem-dsIII-6-8}, to show that for all
  measurable subsets $E\subset H$,
  \begin{equation}\label{eq-b5}
    \phi_{E}\Bigl(\intlohc_{K}g(t)\, dt\Bigr) = \phi_{E} \Bigl(s\mapsto
    \intc_{K} f(t,s)\, dt\Bigr).
  \end{equation}
  Since bounded linear maps commute with integrals, the left-hand side
  of \eqref{eq-b5} is
  \begin{align*}
    \intc_{K}\phi_{E}\bigl(g(t)\bigr)\, dt &= \intc_{K}\intc_{E}
    g(t)(s) \, ds\, dt \\
    \intertext{which, by assumption, is}
    &=\intc_{K}\intc_{E}f(t,s) \, ds\, dt \\
    \intertext{which, since having $f\in\L^{1}(K\times E,C)$ allows us
      to apply Fubini's Theorem, is}
    &= \intc_{E}\intc_{K} f(t,s)\, dt\, ds \\
    &= \phi_{E}\Bigl( s\mapsto \intc_{K}f(t,s)\, dt\Bigr).
  \end{align*}
  This establishes \eqref{eq-b5} and completes the proof.
\end{proof}

\begin{example}
  \label{ex-convolution}
  Now suppose that $(C,H,\alpha)$ is a dynamical system, and that $h$
  and $k$ are in $\slohc$.  Then
  \begin{equation*}
    h*k=\intlohc_{H} h(r)i_{H}(r)(k)\, dr.
  \end{equation*}
\end{example}

\begin{proof}[Proof of the Example]
  We want to apply Lemma~\ref{lem-dunford} with
  \begin{equation*}
    g(r):= h(r)i_{H}(r)(k)\quad\text{and}\quad f(r,s):= h(r)
    i_{H}(r)(k)(s).
  \end{equation*}
  In order to do so, we have to check that $g$ and $f$ are integrable.
  However, $f(r,s) = h(r)\alpha_{r}\bigl(k(r^{-1}s)\bigr)$, and $f$ is
  known to be in $\L^{1}(K\times H,C)$ by standard arguments.
  
  Note that $r\mapsto i_{H}(r)(k)$ is continuous from $H$ to $\lohc$.
  It follows that $r\mapsto h(r) i_{H}(r)(k)$ is measurable from $H$
  to $\lohc$ and has $\sigma$-finite support.  Since
  $\|h(r)i_{H}(r)(k)\|\le\|k\|_{1}\|h(r)\|$, it follows from Tonelli
  that $g\in \L^{1}\bigl(K,\lohc\bigr)$.
  
  Now the result follows immediately from Lemma~\ref{lem-dunford}: for
  almost all $s$,
  \begin{equation*}
    \intlohc_{H}g(r)\, dr(s)=\intc_{H} f(r,s)\, dr = h*k(s). \qed
  \end{equation*}
  \renewcommand{\qed}{}
\end{proof}

\begin{remark}
  \label{rem-notation}
  In the statement of Lemma~\ref{lem-variable}, we are viewing the
  $e_{i}$ as elements of $C\rtimes_{\tau,r}H$, and consequently we
  write their product as $e_{1}e_{2}$.  In the proof however, we will
  want to use that each $e_{i}$ is in $L^{1}(H,C)$, and
  that the product is given by convolution.  Thus it will be a bit
  clearer to use the usual notation $e_{1}*e_{2}$ for their product.
\end{remark}
\begin{proof}
  [Proof of Lemma~\ref{lem-variable}] Let
  \begin{equation}
    \label{eq-b4}
    h(r,t,s):= e_{1}(r)i_{H}(r)\bigl(\sigmat_{t}(e_{2})\bigr)(s).
  \end{equation}
  Our assumptions imply that $h$ is integrable on $H\times H\times K$:
  for example if $e_{1}(r) = \Delta_{H}(r)^{-\frac12}
  u\tau_{r}(v^{*})$ and
  $e_{2}(r)=\Delta_{H}(r)^{-\frac12}w^{*}\tau_{r}(x)$, then
\begin{align*}
  h(r,t,s)&= \Delta_{H}(r)^{-\frac12} u\tau_{r}(v^{*}) \tau_{r}\bigl(
  \sigmat_{t}(e_{2})\bigr)(r^{-1}s) \\
  &= \Delta_{H}(r)^{-\frac12} u\tau_{r}(v^{*})
  \tau_{r}\bigl(\sigma_{t}\bigl(
  w^{*}\tau_{r^{-1}s}(x)\bigr)\bigr)\Delta_{H}(r^{-1}s)^{-\frac12} \\
  &= u\tau_{r}(v^{*})
  \sigma_{t}\bigl(\tau_{r}(w^{*})\tau_{s}(x)\bigr)\Delta_{H}(s)^{-\frac12},
\end{align*}
which is integrable by assumption.  In general, the integrability of
$h$ follows as $e_{1}$ and $e_{2}$ are sums of functions of the form
given above.

Now Fubini's Theorem implies that $r\mapsto h(r,s,t)$ is integrable
for almost all $(t,s)$, and that
\begin{equation}\label{eq-b9}
  f(t,s):=\intc_{H} h(r,t,s)\, dr
\end{equation}
defines an integrable function $f$ on $K\times H$.\footnote{Of course,
  \eqref{eq-b9} only defines $f$ almost everywhere.  As is standard
  practice, we assume the convention that $f(t,s):=0$ when $r\mapsto
  h(r,t,s)$ is not integrable.}  Furthermore, it follows from
Example~\ref{ex-convolution} that for almost all $t$ and $s$,
\begin{equation*}
  f(t,s)=e_{1}*\sigmat_{t}(e_{2})(s).
\end{equation*}

Next we define $g(t):= e_{1}*\sigmat_{t}(e_{2})$.  Then $g$ is a
function from $K$ to $\slohc$ and $g(t)=f(t,\cdot)$ for almost all
$t$.  We want to see that $g\in \L^{1}\bigl(K,\lohc\bigr)$.  However,
$t\mapsto\sigmat_{t}(e_{2})$ is continuous from $K$ to $\lohc$.
  %\footnote{We need to show that $t\mapsto \sigmat_{t}(f)$ is
  % continuous for each $f\in\lohc$.  But it certainly
  % suffices to take $f\in C_{c}(H,C)$, and then it is easy to see that
  % $\sigmat_{t}(f)\to f$ as $t\to e$ in $K$.}  
Therefore $g$ itself is continuous, hence measurable, and it suffices
to see that $\|g\|$ is integrable.  But
\begin{align*}
  \int_{K}\|g(t)\|_{1}\, dt &=
  \int_{K}\int_{H}\|g(t)(s)\|\, ds\, dt \\
  &= \int_{K}\int_{H}\|f(t,s)\|\, ds\, dt \\
  &= \|f\|_{1}<\infty.
\end{align*}
Therefore $g\in \L^{1}\bigl(K,\lohc\bigr)$.  Now an application of
Lemma~\ref{lem-dunford} implies that \eqref{eq:8} is a representative
for
\begin{equation*}
  \intlohc_{K} g(t)\, dt.
\end{equation*}
Since inclusion is a bounded linear map of $\lohc$ into
$C\rtimes_{\tau,r}H$, it follows that \eqref{eq:8} is also a
representative for
\begin{equation*}
  \int^{C\rtimes_{\tau,r}H}_{K}g(t)\, dt,
\end{equation*}
which is what we wanted to prove.
\end{proof}

\section{The proof of Lemma~\ref{lem-dana-sat}}\label{app3}

\begin{lemma}[Rieffel]
  Suppose $G$ acts freely and properly on a locally compact Hausdorff
  space $P$, there is a nondegenerate homomorphism $\theta:C_0(P)\to
  M(A)$, and an action $\alpha$ of $G$ on $A$ such that
  $\alpha_s(\theta(f)a)=\theta(\lt_s(f))\alpha_s(a)$.  Then the proper
  action $\alpha$ of $G$ on $A$ is saturated with respect to
%  $A_1:=\theta(C_c(P))A_0\theta(C_c(P))$, where
  $A_0:=\theta(C_c(P))A\theta(C_c(P))$.
\end{lemma}

  For the proof, we need to know that the action of $G$ on $C_{0}(P)$
  is saturated in the sense of \cite{rie:pm88}; this is proved, for
  example, in \cite{rie:pspm82}.
  
  Let $E_{0}=\lip E<A_{0}, A_{0}>\subset L^1(G,A)$ where $E$ is the
  closure of $E_{0}$ in $A\rtimes_{\alpha,r}G$.  It suffices to see
  that $E_{0}$ is dense in $L^1(G,A)$ in the inductive limit topology.
  To see this, it suffices to see that functions of the form $s\mapsto
  (\phi\otimes a)(s):=\phi(s)a$, where $\phi\in C_c(G)$ and $a$ is in
  any dense subset of $A$, belong to $E$.  Since $A_{0}^{2}$ is dense
  in $A$, we can fix $f,p,k\in C_{c}(P)$, $a,b\in A_0$ and $\phi\in
  C_{c}(G)$ and let
\begin{equation}
  \label{eq:1}
  F(s):= \phi(s)\theta(f)a\theta(p)b^{*}\theta(\bar k),
\end{equation}
and it will suffice to see that $F$ can be approximated in the
inductive limit topology by elements of $E_{0}$.

First we prove a related statement: if $\phi$, $a$, $b$, $f$, $p$ and
$k$ are as above then
  \begin{equation*}
    s\mapsto
    \phi(s)\theta(f)a\theta(p)\alpha_{s}(b^{*})\theta\bigl(\lt_{s}(\bar
    k)\bigr)
  \end{equation*}
  belongs to $E$.  Let $W$ be a compact neighborhood of $\supp \phi$.
  Because the action of $G$ on $C_0(P)$ is proper and saturated with
  respect to $C_c(P)$, given $\epsilon>0$, we can find $g_{i},h_{i}\in
  C_{c}(P)$ such that
  \begin{gather*}
    \label{eq:2}
    \Bigl\| \phi(s) p- \Delta(s)^{-1/2}\sum_{i} g_{i}\lt_{s}(\bar
    h_{i}) \Bigr\|_{\infty}< \frac{\epsilon}
    {\|f\|_{\infty}\|a\|\|b\|\|k\|_{\infty}}
    \quad\text{and} \\
    \supp \bigl(s\mapsto \Delta(s)^{-1/2}\sum_{i} g_{i}\lt_{s}(\bar
    h_{i})\bigr) \subset W.
  \end{gather*}
  Now let
\begin{equation*}
  a_{i}:= \theta(f)a\theta(g_{i})\quad\text{and}\quad b_{i}:=
  \theta(k)b\theta(h_{i}). 
\end{equation*}
Then
\begin{align*}
  \sum_{i}\lip E<a_{i}, b_{i}>(s) &= \sum_{i}\Delta(s)^{-1/2}a_{i}
  \alpha_{s}(b_{i}^{*}) \\
  &=\theta(f)a\theta\Bigl(\Delta(s)^{-1/2} \sum_{i} g_{i}\lt_{s}(\bar
  h_{i}) \Bigr) \alpha_{s}(b^{*})\theta\bigl(\lt_{s}(\bar k)\bigr).
\end{align*}
Now
\begin{equation*}
  \Bigl\| \theta\Bigl(\Delta(s)^{-\half}\sum_{i} g_{i}\lt_{s}(\bar
  h_{i}) \Bigr) - \phi(s)\theta(p)\Bigr\| <
  \frac\epsilon{\|f\|_{\infty}\|a\| \|b\|
  \|k\|_{\infty}}.
\end{equation*}
Thus
\begin{equation*}
  \Bigl\| \sum_{i}\lip E<a_{i},b_{i}> (s) - \phi(s)\theta(f)a\theta(p)
  \alpha_{s} (b^{*}) \theta\bigl(\lt_{s}(\bar k)\bigr) \Bigr\| <\epsilon.
\end{equation*}
Since the neighborhood $W$ does not depend on $\epsilon$, and
$\epsilon$ is arbitrary, it follows that the function
$s\mapsto\phi(s)\theta(f)a\theta(p)\alpha_{s}(b^{*})\theta\bigl(\lt_{s}(\bar
k)\bigr)$ is in $E$.

Now let $N$ be a neighborhood of $e$ in $G$ such that $s\in N$ implies
that
  \begin{equation*}
    \|b^{*}\theta(\bar k) - \alpha_{s}(b^{*})\theta\bigl(\lt_{s}( \bar
    k)\bigr) \| = \|b^{*}\theta(\bar k) -
    \alpha_{s}\bigl(b^{*}\theta(\bar k)\bigr)\|<
    \frac\epsilon{\|\phi\|_{\infty}\|f\|_{\infty} \|a\| \|p\|_{\infty} }.
  \end{equation*}
  Choose $r_{1},\dots,r_{n}\in G$ such that $\supp\phi\subset \bigcup
  Nr_{i}$.  Let $\{\phi_{i}\}\subset C_{c}^{+}(G)$ be such that $\supp
  \phi_{i}\subset Nr_{i}$ and $\sum_{i}\phi_{i}\equiv 1 $ on $\supp
  \phi$ and dominated by $1$ elsewhere.  We showed above that
\begin{equation*}
  F_{i}(s):=
  \phi(s)\phi_{i}(s)\theta(f)a\theta(p)\alpha_{sr_{i}^{-1}}(b^{*})
  \theta\bigl( \lt_{sr_{i}^{-1}}(\bar k)\bigr)
\end{equation*}
defines an element of $E$; if
\begin{equation*}
  F(s):= \phi(s)\theta(f)a\theta(p)b^{*}\theta(\bar k),
\end{equation*}
then
\begin{align*}
  \Bigl\|F(s)-\sum F_{i}(s)\Bigr\| & =
  \Bigl\| \sum \phi_{i}(s)F(s) - F_{i}(s)\Bigr\| \\
  &= \Bigl\| \sum \phi(s)\phi_{i}(s)\theta(f) a \theta(p) \bigl(
  b^{*}\theta(\bar k) - \alpha_{sr_{i}^{-1}}(b^{*}\bar k)\bigr)
  \Bigr\|
  \\
  &\le \sum
  \phi_{i}(s)\|\phi\|_{\infty}\|f\|_{\infty}\|a\|\|p\|_{\infty}
  \|b^{*} \theta(\bar k) - \alpha_{sr_{i}^{-1}}\bigl(b^{*}\theta(\bar
  k)\bigr) \| \\
  \intertext{which, since we may assume $s\in Nr_{i}$, is}
  &\le\epsilon\sum \phi_{i}(s)\le\epsilon.
\end{align*}
Since $\sum_{i}F_{i}\in E$ and both $\supp F_{i}$ and $ \supp F\subset
\supp\phi$,  the result follows.

%%%%%%%%%%%%%%%%%%%%%

%===================================================================
%   REFERENCES
%=================================================================
%   BiBTeX Version
%================================================================
%\bibliographystyle{amsalpha}      % For drafts
% \bibliographystyle{amsplain}      % For final copy
% \bibliography{refdb}
%=================================================================
%   By Hand
%=================================================================
\def\noopsort#1{}\def\cprime{$'$}
\providecommand{\bysame}{\leavevmode\hbox to3em{\hrulefill}\thinspace}
\providecommand{\MR}{\relax\ifhmode\unskip\space\fi MR }
\providecommand{\MRhref}[2]{%
  \href{http://www.ams.org/mathscinet-getitem?mr=#1}{#2} }
\providecommand{\href}[2]{#2}

\end{document}